\documentclass[reqno,12pt]{amsart}
\usepackage{epsfig}
\usepackage{amsfonts,amsmath,amssymb,mathrsfs,verbatim,amsthm,amscd}
\usepackage{cite}
\def\II{\hbox{{1}\kern-.25em\hbox{l}}}
\newcommand{\V}{\mathbb V}
\newcommand{\be}{\begin{equation}}
\newcommand{\ee}{\end{equation}}
\newcommand{\ba}{\begin{eqnarray}}
\newcommand{\ea}{\end{eqnarray}}

\newcommand{\R}{\mathbb R}
\newcommand{\Z}{\mathbb Z}

\begin{document}

\title[Elliptic modular double]
{Finite-dimensional representations \\ of the elliptic modular double}

\author{ S. E. Derkachov}
 \address{PDMI RAS, Fontanka 27, St. Petersburg, Russia}

\author{V. P. Spiridonov}
 \address{Laboratory of Theoretical Physics, JINR, Dubna, Moscow region, 141980, Russia}

\begin{abstract}
We investigate the kernel space of an integral operator $\mathrm{M}(g)$
depending on the ``spin" $g$ and describing an elliptic Fourier transformation.
The operator $\mathrm{M}(g)$ is an intertwiner for the elliptic modular double
formed from a pair of Sklyanin algebras with the parameters $\eta$ and
$\tau$, Im$\,\tau>0$, Im$\,\eta>0$. For two-dimensional lattices
$g=n\eta + m\tau/2$ and  $g=1/2+n\eta + m\tau/2$
with incommensurate $1, 2\eta,\tau$ and integers $n,m>0$,
the operator $\mathrm{M}(g)$
has a finite-dimensional kernel that consists of the products of theta functions
with two different modular parameters and is invariant under the action
of generators of the elliptic modular double.
\end{abstract}

\maketitle


\medskip

\hfill{\em To Ludwig Faddeev in honor of his 80th birthday}

\medskip

\tableofcontents

\section{An elliptic modular double}

The concept of a modular double of a quantum group was introduced
by Faddeev in \cite{fad:mod}, where it was shown that the
quantum algebra $U_q(sl_2)$ with the deformation parameter $q$
does not uniquely define its representation space and requires an extension.
Such an extension is formed by two sets of generators
$\mathbf{E}\,,\mathbf{F}\,,\mathbf{K}$ and
$\widetilde{\mathbf{E}}\,,\widetilde{\mathbf{F}}\,,\widetilde{\mathbf{K}}$.
The usual algebraic relations
\begin{equation} \label{qsl2}
\begin{array}{c}
[\mathbf{E},\mathbf{F}] = \frac{\mathbf{K}^2 - \mathbf{K}^{-2}}{q-q^{-1}} \;,\;\;\;
\mathbf{K} \mathbf{E} = q \mathbf{E} \mathbf{K} \;,\;\;\;
\mathbf{K} \mathbf{F} = q^{-1} \mathbf{F} \mathbf{K},
\end{array}
\end{equation}
where $q = e^{\pi \textup{i}  \sigma}$, are supplemented by similar relations
for the operators $\widetilde{\mathbf{E}},\widetilde{\mathbf{F}},$ and
$\widetilde{\mathbf{K}}$
with the parameter $q$ replaced with $\widetilde{q} = e^{\pi \textup{i}/\sigma}$.
The generators
$\mathbf{E}$ and $\mathbf{F}$ commute with $\widetilde{\mathbf{E}}$ and
$\widetilde{\mathbf{F}}$. The generator
$\mathbf{K}$ anticommutes with $\widetilde{\mathbf{E}}$ and $\widetilde{\mathbf{F}}$,
and $\widetilde{\mathbf{K}}$ anticommutes with $\mathbf{E}$ and $\mathbf{F}$.
This algebra has two central elements, Casimir operators.
The first of them has the form
\begin{equation} \label{Casimir}
\mathbf{C}= \left(q-q^{-1}\right)^2\,\mathbf{F} \mathbf{E} - q \mathbf{K}^2
- q^{-1} \mathbf{K}^{-2} + 2\,.
\end{equation}
The second is constructed from
$\widetilde{\mathbf{E}},\widetilde{\mathbf{F}},\widetilde{\mathbf{K}}$,
and $\widetilde{q}$; its explicit expression is similar to (\ref{Casimir}).
Particular representations of this
modular double were considered in \cite{bult,BT,fad:mod,fad:mod2,KLS,PT}.

Quantum algebras emerged from the theory of the Yang-Baxter equation (YBE)
\begin{equation}\label{YBE}
\mathbb{R}_{12} (u-v)\,\mathbb{R}_{13}(u)\, \mathbb{R}_{23}(v)
=\mathbb{R}_{23}(v)\,\mathbb{R}_{13}(u)\,\mathbb{R}_{12}(u-v),
\end{equation}
which plays an important role in mathematical physics~\cite{Baxter,FT,Jimbo}.
Here, the operators $\mathbb{R}_{jk}$ act in the subspace $\V_j\otimes\V_k$
of the tensor product $\V_1\otimes\V_2\otimes\V_3$ of three (in general different)
spaces $\V_j$. Variables $u$ and $v$ are called spectral parameters.
The Sklyanin algebra \cite{skl1,skl2} is a one parameter deformation of $U_q(sl_2)$
or an elliptic deformation of the $sl_2$-algebra. It emerges
from equation \eqref{YBE} when $\mathbb{R}_{12} (u)$ is given by Baxter's
$4\times 4$ R-matrix \cite{Baxter},
\begin{equation}\label{Baxter}
\mathbb{R}_{12}(u) = \sum_{a=0}^3 w_{a} (u)\,
\sigma_a \otimes\sigma_a \ \ ,\ \ w_{a}(u) = \frac{\theta_{a+1}
(u+\eta|\tau)}{\theta_{a+1}(\eta|\tau)}\,,
\end{equation}
and $\mathbb{R}_{13}(u)$, $\mathbb{R}_{23}(v)$
are  $2\times 2$ matrices fixed as copies of the L-operator:
\begin{equation}\makebox[-1em]{}
\mathrm{L}(u):=
\sum_{a=0}^3 w_{a} (u)\, \sigma_a \otimes
\mathbf{S}^a = \left(
\begin{array}{cc}
w_0(u)\,\mathbf{S}^0+w_3(u)\,\mathbf{S}^3 &
w_1(u)\,\mathbf{S}^1-\textup{i} w_2(u)\,\mathbf{S}^2 \\
w_1(u)\,\mathbf{S}^1+\textup{i} w_2(u)\,\mathbf{S}^2&
w_0(u)\,\mathbf{S}^0-w_3(u)\,\mathbf{S}^3
\end{array} \right),
\label{L_op}\end{equation}
where $\sigma_a$ are the Pauli matrices and $\eta$ and $\tau$ are free
parameters, Im$\,\tau>0$. Here, $\theta_{a}(u|\tau)$ are the Jacobi theta-functions
\begin{eqnarray}\label{theta1} &&
\theta_{1}(z|\tau) = -\sum_{n\in\mathbb{Z}}
\mathrm{e}^{\pi \textup{i} \left(n+\frac{1}{2}\right)^2\tau}\cdot
\mathrm{e}^{2\pi \textup{i}
\left(n+\frac{1}{2}\right)\left(z+\frac{1}{2}\right)},
\qquad \theta_{2}(z|\tau)=\theta_1(z+{\textstyle\frac{1}{2}}|\tau),
\\ &&
\theta_{3}(z|\tau)=e^{\frac{\pi \textup{i}\tau}{4}+\pi \textup{i} z}
\theta_2(z+{\textstyle \frac{\tau}{2}}|\tau), \qquad
\theta_4(z|\tau)= \theta_3(z+{\textstyle\frac{1}{2}}|\tau).
\nonumber\end{eqnarray}
It is convenient to introduce the basic variables $p=e^{2\pi \textup{i} \tau}$
and $q=e^{4\pi \textup{i} \eta}$. Then
$$
\theta_{1}(z|\tau) =
\frac{e^{-\pi\textup{i}z} \theta(e^{2\pi \textup{i}z};p)}{\mathrm{R}(\tau)}, \quad
\mathrm{R}(\tau) = \frac{p^{-\frac{1}{8}}}{\textup{i} (p;p)_\infty},
$$
where
$$
\theta(t;p)=(t;p)_\infty(pt^{-1};p)_\infty,\qquad
(t;p)_\infty=\prod_{k=0}^\infty(1-tp^k).
$$
Theta functions are quasiperiodic, for example
$$
\theta_{1}(z+1|\tau) =-\theta_{1}(z|\tau), \quad
\theta_{1}(z+\tau|\tau) =- e^{-2\pi\textup{i} z-\pi\textup{i}\tau}\theta_{1}(z|\tau).
$$
In what follows, we need the general relations
\begin{equation}\label{mper}
\theta_a(z+ m\,\tau\,|\tau) =
\mu_a\,
e^{-\pi \textup{i}\tau\,m^2- 2\pi \textup{i} \,m z}\,\theta_a(z\,|\tau)\,,
\quad m\in\Z,
\end{equation}
where $\mu_a=(-1)^m$ for $a=1,4$ and $\mu_a=1$ for $a=2,3$.

The YBE corresponding to operators \eqref{Baxter} and \eqref{L_op}
takes the form of the RLL-relation
$$\mathbb{R}_{12} (u-v)\,\mathrm{L}_{1}(u)\, \mathrm{L}_{2}(v)
=\mathrm{L}_{2}(v)\,\mathrm{L}_{1}(u)\,\mathbb{R}_{12}(u-v)\, ,
$$
which yields the following Sklyanin algebra \cite{skl1}:
\begin{eqnarray}\nonumber &&
\mathbf{S}^\alpha\,\mathbf{S}^\beta - \mathbf{S}^\beta\,\mathbf{S}^\alpha =
\textup{i}\left(\mathbf{S}^0\,\mathbf{S}^\gamma +\mathbf{S}^\gamma\,\mathbf{S}^0\right)\,,
\\ &&
\mathbf{S}^0\,\mathbf{S}^\alpha - \mathbf{S}^\alpha\,\mathbf{S}^0 =
\textup{i}\,\mathbf{J}_{\beta \gamma}\left(\mathbf{S}^\beta\,\mathbf{S}^\gamma +\mathbf{S}^\gamma\,\mathbf{S}^\beta\right)\,,
\label{sklalg}\end{eqnarray}
where the triplet $(\alpha,\beta,\gamma)$ is an arbitrary cyclic permutation
of $(1,2,3)$. The structure constants $\mathbf{J}_{\alpha\beta}$ are not independent
and satisfy  the constraint
\begin{equation}
\mathbf{J}_{12}+\mathbf{J}_{23}+\mathbf{J}_{31}
+\mathbf{J}_{12}\mathbf{J}_{23}\mathbf{J}_{31}=0.
\label{con}\end{equation}
Parametrizing them as $\mathbf{J}_{\alpha\beta}
=\frac{\mathbf{J}_{\beta}-\mathbf{J}_{\alpha}}{\mathbf{J}_{\gamma}}$ for
$\gamma\neq \alpha,\beta$, which automatically resolves the constraint \eqref{con},
we can write
\begin{equation}
\mathbf{J}_{1}=\frac{\theta_2(2\eta|\tau)\theta_2(0|\tau)}
{\theta_2^2(\eta|\tau)},\quad
\mathbf{J}_{2}=\frac{\theta_3(2\eta|\tau)\theta_3(0|\tau)}
{\theta_3^2(\eta|\tau)},\quad
\mathbf{J}_{3}= \frac{\theta_4(2\eta|\tau)\theta_4(0|\tau)}
{\theta_4^2(\eta|\tau)}.
\label{uni1}\end{equation}
The constants $\mathbf{J}_{\alpha}$ are also not independent, because they
are parametrized by only two complex variables $\eta$ and $\tau$.

There are two Casimir operators commuting with all generators:
$$
\mathbf{K}_0 = \sum_{a=0}^3\,\mathbf{S}^a\,\mathbf{S}^a\ ,\qquad
\mathbf{K}_2 = \sum_{\alpha=1}^3\,\mathbf{J}_\alpha\,
\mathbf{S}^\alpha\,\mathbf{S}^\alpha\, ,
$$
$$
\left[\mathbf{K}_0 ,\mathbf{S}^a\right] =
\left[\mathbf{K}_2 ,\mathbf{S}^a\right] =0.
$$

The operators $\mathbf{S}^a$ can be realized as finite-difference operators
acting on functions of a complex variable $z$ \cite{skl2}:
\begin{eqnarray}\nonumber &&\makebox[-3em]{}
\mathbf{S}^a= e^{\pi\textup{i}z^2/\eta}\frac{\textup{i}^{\delta_{a,2}}
\theta_{a+1}(\eta|\tau)}{\theta_1(2 z|\tau) } \Bigl[\,\theta_{a+1} \left(2
z-g +\eta|\tau\right)e^{\eta\partial_z} - \theta_{a+1}
\left(-2z-g+\eta|\tau\right)e^{-\eta\partial_z}\Bigl]e^{-\pi\textup{i}z^2/\eta},
\\  && \makebox[0em]{}
=f_a(z)e^{\eta\partial_z}+f_a(-z)e^{-\eta\partial_z},\quad
\label{Sklyan} \end{eqnarray}
where
$$
f_a(z)=e^{-\pi\textup{i}\eta-2\pi\textup{i}z}\, \textup{i}^{\delta_{a,2}}
\theta_{a+1}(\eta|\tau)\frac{\theta_{a+1} \left(2
z-g +\eta|\tau\right)}{\theta_1(2 z|\tau)}
$$
and $e^{\pm\eta\partial_z}$ denote the shift operators,
$e^{\pm\eta\partial_z}f(z)=f(z\pm\eta)$. The variable $g$ is usually defined
as $g=\eta(2\ell+1)$ and the parameter $\ell\in\mathbb{C}$
is called the spin. We also call $g$ the spin.
The parameters $\tau$, $\eta$, and $g$ characterize
representations of the Sklyanin algebra because they fix the
values of the Casimir operators. We note that our
operators \eqref{Sklyan} differ from the standard ones by
multiplication by exponentials $e^{\pm\pi\textup{i}z^2/\eta}$ from the left
and right (such a choice leads to the analyticity of the intertwining operator
in $e^{2\pi\textup{i} z}$; see \cite{DS} and the considerations below).

The {\em elliptic modular double} was introduced in \cite{AA2008}.
A particular realization of this algebra degenerates
to Faddeev's modular double  in a special limit \cite{fad:mod}.
The double of interest for us is generated by $\mathbf{S}^a$ and
new operators $\mathbf{\tilde S}^a$ satisfying the relations
\begin{eqnarray}\nonumber &&
\mathbf{\tilde S}^\alpha\,\mathbf{\tilde S}^\beta - \mathbf{\tilde S}^\beta\,\mathbf{\tilde S}^\alpha =
\textup{i}\left(\mathbf{\tilde S}^0\,\mathbf{\tilde S}^\gamma +\mathbf{\tilde S}^\gamma\,\mathbf{\tilde S}^0\right)\,,
\\ &&
\mathbf{\tilde S}^0\,\mathbf{\tilde S}^\alpha - \mathbf{\tilde S}^\alpha\,\mathbf{\tilde S}^0 =
\textup{i}\,\mathbf{\tilde J}_{\beta \gamma}\left(\mathbf{S}^\beta\,\mathbf{S}^\gamma +\mathbf{S}^\gamma\,\mathbf{S}^\beta\right)\,,
\label{sklalg_doub}\end{eqnarray}
where the triplet $(\alpha,\beta,\gamma)$ is an arbitrary cyclic permutation
of $(1,2,3)$. The tilded structure constants have a similar parametrization
$\mathbf{\tilde J}_{\alpha\beta}=\frac{\mathbf{\tilde J}_{\beta}-\mathbf{\tilde J}_{\alpha}}
{\mathbf{\tilde J}_{\gamma}}$ for $\gamma\neq \alpha,\beta$,  where
\begin{equation}
\mathbf{\tilde J}_{1}=\frac{\theta_2(\tau|2\eta)\theta_2(0|2\eta)}
{\theta_2^2(\tau/2|2\eta)},\quad
\mathbf{\tilde J}_{2}=\frac{\theta_3(\tau|2\eta)\theta_3(0|2\eta)}
{\theta_3^2(\tau/2|2\eta)},\quad
\mathbf{\tilde J}_{3}= \frac{\theta_4(\tau|2\eta)\theta_4(0|2\eta)}
{\theta_4^2(\tau/2|2\eta)}.
\label{uni2}\end{equation}
We note that the constants $\mathbf{\tilde J}_{\alpha}$ depend on
$\mathbf{J}_{\alpha}$. Uniformizations
\eqref{uni1} and \eqref{uni2} give an implicit form of this dependence;
it would be interesting to describe it purely algebraically or geometrically.

There are two Casimir operators commuting with all generators:
$$
\mathbf{\tilde K}_0 = \sum_{a=0}^3\,\mathbf{\tilde S}^a\,\mathbf{\tilde S}^a\ ,\qquad
\mathbf{\tilde K}_2 = \sum_{\alpha=1}^3\,\mathbf{\tilde J}_\alpha\,\mathbf{\tilde S}^\alpha\,
\mathbf{\tilde S}^\alpha\,,
$$
$$
\left[\mathbf{\tilde K}_0 ,\mathbf{\tilde S}^a\right] =
\left[\mathbf{\tilde K}_2 ,\mathbf{\tilde S}^a\right] =0.
$$
The cross-commutation relations between $\mathbf{S}^a$ and $\mathbf{\tilde S}^a$ have
the forms
\begin{eqnarray} \nonumber
&& \mathbf{S}^a \mathbf{\tilde S}^b=\mathbf{\tilde S}^b \mathbf{S}^a, \quad a,b\in\{0,3\} \quad
\text{or} \quad a,b\in\{1,2\},
\\
&& \mathbf{S}^a \mathbf{\tilde S}^b=-\mathbf{\tilde S}^b \mathbf{S}^a, \quad a\in\{0,3\},\; b\in\{1,2\} \quad
\text{or} \quad a\in\{1,2\},\;b\in\{0,3\}.
\end{eqnarray}
It follows that $\mathbf{K}_0$ and $\mathbf{K}_2$ commute with
$ \mathbf{\tilde S}^b$ and, vice versa,
 $\mathbf{\tilde K}_0$ and $\mathbf{\tilde K}_2$ commute with $ \mathbf{S}^b$, i.e., we
have four  Casimir operators:
$$
[\mathbf{K}_0,\mathbf{\tilde S}^a]=[\mathbf{K}_2,\mathbf{\tilde S}^a]=
[\mathbf{\tilde K}_0, \mathbf{S}^a]=[\mathbf{\tilde K}_2, \mathbf{S}^a]=0.
$$
The elliptic modular double is an associative algebra, similar to the standard
Sklyanin algebra. Because the generators $\mathbf{S}^a$
and $\mathbf{\tilde S}^a$ do not commute, it is not a direct product of two Sklyanin
algebras, although it is quite easy to trace the difference
of actions of the subalgebra generators on modules in different orders.

There is an obvious finite-difference operator
realization of the second Sklyanin algebra generators.
Namely, the operators $\mathbf{\tilde S}^a$ are obtained
from $\mathbf{S}^a$ by simply permuting $2\eta$ and $\tau$
(which requires Im$\, \eta>0$):
\begin{eqnarray}\nonumber &&
\mathbf{\tilde S}^a  = e^{2\pi\textup{i}z^2/\tau}
\frac{\textup{i}^{\delta_{a,2}}\theta_{a+1}(\tau/2|2\eta)}{\theta_1(2 z|2\eta) }
 \Bigl[\,\theta_{a+1} \left(2z-g+\frac{\tau}{2}\big|2\eta\right)
 \mathrm{e}^{\frac{1}{2}\tau \partial_z}
\\ && \makebox[10em]{}
 - \theta_{a+1}\left(-2z-g+\frac{\tau}{2}\big|2\eta\right)
\mathrm{e}^{-\frac{1}{2}\tau \partial_z}\, \Bigl]e^{-2\pi\textup{i}z^2/\tau},
\label{mod_doub2}\end{eqnarray}
where $g$ is the same arbitrary parameter as in \eqref{Sklyan}.

In the considered realization, the Casimir operators reduce to the
scalar expressions
\begin{eqnarray*} &&
\mathbf{K}_0 = 4\,\theta_1^2\bigl(g|\tau\bigr)\ ,\quad
\mathbf{K}_2 = 4\,\theta_1\bigl(g-\eta|\tau\bigr)\,
\theta_1(g+\eta|\tau)\,,
\\ &&
\mathbf{\tilde K}_0 = 4\,\theta_1^2\bigl(g|2\eta\bigr)\ ,\quad
\mathbf{\tilde K}_2 = 4\,\theta_1\bigl(g-{\textstyle \frac{\tau}{2}}|2\eta\bigr)\,
\theta_1(g+{\textstyle \frac{\tau}{2}}|2\eta)\,,
\end{eqnarray*}
which  are invariant under the reflection $g\to -g$.
The variables $\eta$ and $\tau$ are fixed by the structure constants,
and the spin parameter $g$ fixes the values of
all Casimir operators and determines representations of the
elliptic modular double. There exists another elliptic modular double
\cite{AA2008} based on the transformation
$\eta \to 1/(4\eta),\, \tau \to \tau/(2\eta), \, z \to z/(2\eta), \,
g \to g/(2\eta)$ in definition \eqref{Sklyan}, for which the
regime Im$\,\eta=0$ is permitted, but we do not consider it here.

\section{An intertwining operator}

Here we focus our attention on the integral operator introduced in \cite{spi:bailey}
for defining a universal integral transform of hypergeometric type
yielding an integral analogue of the Bailey chain techniques \cite{aar}.
This operator acts on holomorphic functions of one complex variable $\Phi(x)$ as
\begin{equation}
[\mathrm{M}(g)\Phi](z)= \frac{(q;q)_\infty\,(p;p)_\infty}{2}
\int_0^1 \frac{\Gamma(\pm z\pm x -g)}
{\Gamma(-2g,\pm 2x)}\Phi(x)dx,
\label{S1fin}\end{equation}
where  Im$(-g\pm z)>0$. We use the notation
 $\Gamma(a,b\pm z):=\Gamma(a)\Gamma(b+z)\Gamma(b-z)$, where
\begin{equation}
\Gamma(z):=\Gamma(z|\tau,2\eta):=
\prod_{n,m=0}^{\infty} \frac{1-\mathrm{e}^{-2\pi\textup{i}z}
p^{n+1}q^{m+1}}{1-\mathrm{e}^{2\pi\textup{i}z} p^ nq^m}
\label{egamma}\ee
is the elliptic gamma function defined for $|p|,|q|<1$.
The constraints on the parameter $g$ and argument $z$ can be relaxed
by deforming the integration contour, i.e., by analytical continuation
of the expression given above (in the
situations when this is allowed by the test functions).

The intertwining operator satisfies a very simple inversion relation
resembling the key Fourier transform property \cite{spi-war:inversions}
\begin{equation}
\mathrm{M}(g)\, \mathrm{M}(-g) = \II\,.
\label{inv}\end{equation}
This equality holds for an appropriate space of test functions,
at least for the values of $g$ away from the two discrete lattices
$g = n\eta+m\frac{\tau}{2}$ and $g = \frac{1}{2} +n\eta+m\frac{\tau}{2},$ where
$n,m\in\mathbb{Z}$ \cite{DS}.

As shown in \cite{DS}, operator \eqref{S1fin}, being symmetric in $2\eta$
and $\tau$, satisfies the intertwining relations:
\begin{equation}
\mathrm{M}(g)\,\mathbf{S}^a(g) =
\mathbf{S}^a(-g)\, \mathrm{M}(g), \qquad
\mathrm{M}(g)\,\mathbf{\tilde S}^a(g) =
\mathbf{\tilde S}^a(-g)\,\mathrm{M}(g)\,.
\label{inter1}\end{equation}
Here, we explicitly indicate the $g$-spin dependence of the Sklyanin algebra generators
in order to show that the parameter $g$ simply changes the sign under the action of
$\mathrm{M}$. Because this does not change the Casimir operator
values, the operator $\mathrm{M}(g)$ connects
equivalent representations to each other.
The notation $g=\eta(2\ell+1)$, which assumes the transformation $\ell\to -1-\ell$,
is used in the conventional Sklyanin algebra setting.
Equalities \eqref{inter1} show that the kernel of the $\mathrm{M}$ -operator forms
an invariant space for the elliptic modular double, i.e., it is invariant under
the action of the Sklyanin algebra generators ${\mathbf S}^a(g)$ \eqref{Sklyan}
and  $\mathbf{\tilde S}^a(g)$ \eqref{mod_doub2}.

We previously discovered a nontrivial finite-dimensional kernel of the operator
$\mathrm{M}$ \cite{DS}. Here, we describe this space in more detail.
The intertwining operators are very important tools
in the representations theory \cite{Gelfand,Kna},
where they are used to analyze reducibility questions like
the existence and characterization of invariant subspaces for a given
infinite-dimensional representation. Also, in the theory of the YBE
\cite{DM}, they can serve as elementary building blocks for the R operators.

The intertwining operator $\mathrm{M}(g)$ plays a key role in building
the most complicated known solution of the YBE \cite{DS} found
along the lines of general construction in \cite{SD,DKK}. The latter approach to the YBE
is based on a twisted representation of the generators of the permutation group.
In our case, the needed Coxeter relations are satisfied as a consequence
of the elliptic beta integral \cite{spi:umn} and the Bailey lemma in
\cite{spi:bailey}. We do not use our results to build
solutions of the YBE here, postponing this task to a separate paper.

We set
$$
Z=e^{2\pi \textup{i} z}, \quad X=e^{2\pi \textup{i} x}, \quad t=e^{-2\pi \textup{i} g}, \quad
\Gamma_{\! p,q}(Z):=\Gamma(z|\tau,2\eta)
$$
and assume that the operator $\mathrm{M}(g)$ acts in the space of holomorphic functions
of $X$. We can then rewrite the intertwining operator in the form
\begin{equation}
[\mathrm{M}(g)f](Z)= \frac{(q;q)_\infty\,(p;p)_\infty}{4\pi \textup{i}}
\int_{\mathbb{T}} \frac{\Gamma_{\! p,q}(tZ^{\pm1}X^{\pm1})}
{\Gamma_{\! p,q}(t^2,X^{\pm 2})}f(X)\frac{dX}{X},
\label{S1fin_m}\end{equation}
where $\mathbb{T}$ is the unit circle of positive orientation, $|tZ^{\pm 1}|<1$, and
$$
\Gamma_{\! p,q}(a,b^{\pm k}):=\Gamma_{\! p,q}(a)\Gamma_{p,q}(b^k)\Gamma_{\! p,q}(b^{-k}).
$$
We note that requiring the holomorphicity in $X$ is equivalent to requiring
periodicity for the functions $\Phi(x)$ in \eqref{S1fin},
$\Phi(x+1)=\Phi(x)$, which strongly restricts the space of test functions.
But such a requirement is natural
because after the action of the intertwining operator $\mathrm{M}$, we
always obtain meromorphic functions of the variable $Z$.

\section{Contiguous relations for the intertwining operator}

Contiguous (or recurrence) relations connect special
functions with different parameter values to each other \cite{aar}.
The first contiguous relation for elliptic hypergeometric
integrals was already constructed in \cite{spi:umn}. Such relations
can also be formulated for integral operators, and we want to
do this here for the intertwining operator $\mathrm{M}(g)$.

The recurrence relation of interest for operator \eqref{S1fin}
has the form \cite{CDKKIII}
\begin{equation}\label{rec1}
\frac{e^{\pi \textup{i} \eta}}{\mathrm{R}(\tau)} \cdot
\frac{e^{\pi \textup{i}\frac{z^2}{ \eta}}}{\theta_1(2z)}
\left[ \bar{\theta}_k (z+g+\eta)\, e^{\eta \partial_z} -
\bar{\theta}_k (z-g-\eta)\, e^{-\eta \partial_z}
\right]\, e^{-\pi \textup{i}\frac{z^2}{ \eta}}\,\mathrm{M}(g) =
 \mathrm{M}(g+\eta) \,\bar{\theta}_k (z)\, ,
\end{equation}
where $\bar{\theta}_k (z)=\theta_k(z|{\textstyle\frac{\tau}{2}}),\, k=3,4$.
Here, $\mathrm{R}(\tau)$ is the constant defined in  \eqref{theta1}.
In the right-hand sides of equality \eqref{rec1} and other
expressions below, we use the variable
$z$ and assume that it is an ``internal" variable, i.e., it plays the role of $x$
in the action of integral operator \eqref{S1fin}.
The given relation can be written as the equality
\begin{equation}
\mathrm{A}_k(g)\,\mathrm{M}(g) = \mathrm{M}(g+\eta)\,\bar{\theta}_k (z)\,,
\label{recrelA}\end{equation}
where $\mathrm{A}_k(g)$ is the difference operator
\begin{equation}
\mathrm{A}_k(g) = \frac{e^{\pi \textup{i} \eta}}{\mathrm{R}(\tau)}\,e^{\pi \textup{i}\frac{z^2}{ \eta}}\,\frac{1}{\theta_1(2z)}
\left[ \bar{\theta}_k (z+g+\eta)\, e^{\eta \partial_z} -
\bar{\theta}_k (z-g-\eta)\, e^{-\eta \partial_z}
\right]\, e^{-\pi \textup{i}\frac{z^2}{ \eta}}\,.
\label{Ak}\end{equation}

To prove the operator identity \eqref{rec1}, we note that
it is equivalent to the equation for the integrand of the intertwining operator
\begin{eqnarray} \nonumber &&
e^{\pi \textup{i}\frac{z^2}{ \eta}}\,\frac{1}{\theta_1(2z)}
\left[ \bar{\theta}_k (z+g+\eta)\, e^{\eta \partial_z} -
\bar{\theta}_k (z-g-\eta)\, e^{-\eta \partial_z}
\right]\, e^{-\pi \textup{i}\frac{z^2}{ \eta}}\,
\frac{\Gamma(\pm z\pm x -g)}{\Gamma(-2g,\pm 2x)} =
\\ && \makebox[2em]{}
= \mathrm{R}(\tau)\, e^{-\pi \textup{i} \eta}\cdot
\frac{\Gamma(\pm z\pm x -g-\eta)}{\Gamma(-2g-2\eta,\pm 2x)}
\,\bar{\theta}_k (x)\,.
\label{rec3}\end{eqnarray}
The proof of this relation is based on two key formulas
\begin{eqnarray} &&
\Gamma(z+2 \eta) = \mathrm{R}(\tau) \, e^{\pi \textup{i} z} \, \theta_1 (z)\, \Gamma(z)\, ,
\label{eqgamma}
\\  \label{theta} &&
2\,\theta_1(x\pm y): = 2\,\theta_1(x+y)\,\theta_1(x-y) = \bar\theta_4(x)\,\bar\theta_3(y)
-\bar\theta_4(y)\,\bar\theta_3(x)\,.
\end{eqnarray}
We have the equalities
$$
e^{\pi \textup{i}\frac{z^2}{ \eta}}\, e^{\eta \partial_z}\,
e^{-\pi \textup{i}\frac{z^2}{ \eta}}\,
\Gamma(\pm z\pm x -g)
= \mathrm{R}^2(\tau)\, e^{-3\pi \textup{i} \eta-2\pi \textup{i} g}
\,\theta_1(z-\eta-g \pm x)\,
\Gamma(\pm z\pm x -g-\eta)\,,
$$
$$
e^{\pi \textup{i}\frac{z^2}{ \eta}}\, e^{-\eta \partial_z}\,
e^{-\pi \textup{i}\frac{z^2}{ \eta}}\,
\Gamma(\pm z\pm x -g)
= \mathrm{R}^2(\tau)\, e^{-3\pi \textup{i} \eta-2\pi \textup{i} g}
\,\theta_1(z+\eta+g \pm x)\,
\Gamma(\pm z\pm x -g-\eta)\,,
$$
and therefore
\begin{eqnarray*} &&
e^{\pi \textup{i}\frac{z^2}{ \eta}}\,\frac{1}{\theta_1(2z)}
\left[ \bar{\theta}_k (z+g+\eta)\, e^{\eta \partial_z} -
\bar{\theta}_k (z-g-\eta)\, e^{-\eta \partial_z}
\right]\, e^{-\pi \textup{i}\frac{z^2}{ \eta}}\,
\Gamma(\pm z\pm x -g) =
\\ && \makebox[3em]{}
 =\mathrm{R}^2(\tau)\, e^{-3\pi \textup{i} \eta-2\pi \textup{i} g}\,\cdot
\frac{1}{\theta_1(2z)}\, \Gamma(\pm z\pm x -g-\eta)\,\cdot
\\ && \makebox[4em]{}
\cdot\,\left[
\bar{\theta}_k (z+g+\eta)\,\theta_1(z-\eta-g \pm x)
-\bar{\theta}_k (z-g-\eta)\,\theta_1(z+\eta+g \pm x)
\right] =
\\ && \makebox[3em]{}
= \mathrm{R}^2(\tau)\, e^{-3\pi \textup{i} \eta-2\pi \textup{i} g}\,
\cdot \Gamma(\pm z\pm x -g-\eta)\,\theta_1(-2g-2\eta)\,\bar{\theta}_k(x)\,.
\end{eqnarray*}
In the last line, we use the equality
\begin{eqnarray}\nonumber &&
\bar{\theta}_k (z+g+\eta)\,\theta_1(z-\eta-g \pm x)
-\bar{\theta}_k (z-g-\eta)\,\theta_1(z+\eta+g \pm x)=
\\ && \makebox[4em]{}
 = \theta_1(2z)\,\theta_1(-2g-2\eta)\,\bar{\theta}_k(x)\,,
\label{id2}\end{eqnarray}
which is derived using relation~(\ref{theta}) twice.
At the final step, we apply the equality
$$
\frac{\theta_1(-2g-2\eta)}{\Gamma(-2g)} = \frac{\mathrm{R}^{-1}(\tau)
\, e^{2\pi \textup{i} (g+\eta)}}{\Gamma(-2g-2\eta)}
$$
and obtain identity ~(\ref{rec3}).

We suppose that the parameter $g$ satisfies the constraints ensuring that
inversion relation \eqref{inv} is satisfied. We can then easily obtain another recurrence
relation:
\begin{equation}\label{rec2}
\bar{\theta}_k (z)\, \mathrm{M}(g+\eta)
=  \mathrm{M}(g) \,
\frac{e^{\pi \textup{i} \eta}}{\mathrm{R}(\tau)}
\frac{e^{\pi \textup{i}\frac{z^2}{ \eta}}}{\theta_1(2z)}
\left[ \bar{\theta}_k(z-g) \,e^{\eta \partial_z} -
\bar{\theta}_k(z+g) \, e^{-\eta \partial_z}
\right]\, e^{-\pi \textup{i}\frac{z^2}{ \eta}}\,.
\end{equation}
Indeed, multiplying equality \eqref{recrelA} by the operator $\mathrm{M}^{-1}(g+\eta)$,
we obtain
$$
\mathrm{M}^{-1}(g+\eta)\,\mathrm{A}_k(g) =
\bar{\theta}_k (z)\,\mathrm{M}^{-1}(g)
\ \ \text{or} \ \
\mathrm{M}(-g-\eta)\,\mathrm{A}_k(g) =
\bar{\theta}_k (z)\,\mathrm{M}(-g).
$$
Changing the sign of $g$ and shifting $g$ by $\eta$ we obtain
\begin{equation}
\mathrm{M}(g)\,\mathrm{A}_k(-g-\eta) =
\bar{\theta}_k (z)\,\mathrm{M}(g+\eta),
\label{RR2}\end{equation}
which coincides with~(\ref{rec2}).
But relation \eqref{rec2} can be derived without using the inversion
relation under some mild restrictions on the test function space.
Indeed, operator relation \eqref{RR2} acting on a function $\Phi(x)$
has the explicit form
\begin{eqnarray}\nonumber
&& \int_0^1\frac{\Gamma(\pm z\pm x -g)}{\mathrm{R}(\tau)\theta_1(2x)\Gamma(\pm 2x)}
\left(\bar \theta_k(x-g) e^{-2\pi \textup{i} x}\Phi(x+\eta)
-\bar \theta_k(x+g) e^{2\pi \textup{i} x}\Phi(x-\eta)\right) dx =
\\ && \makebox[2em]{}
=\frac{\Gamma(-2g)}{\Gamma(-2g-2\eta)}\bar \theta_k(z)\int_0^1
\frac{\Gamma(\pm z\pm x -g-\eta)}{\Gamma(\pm 2x)}\Phi(x)dx.
\label{dualeqn}\end{eqnarray}

We now suppose that the function $\Phi(x)$ is periodic, $\Phi(x+1)=\Phi(x)$,
and that the integrands in the left-hand side of \eqref{dualeqn}
have no simple poles in the parallelogram $[-\eta, \eta, 1+\eta, 1-\eta]$.
This condition means that the integrands
in the left-hand side of \eqref{dualeqn} should be meromorphic functions of
$X=e^{2\pi \textup{i}x}$ without simple poles in the annulus
$ |q|^{1/2} \leq |X|\leq |q|^{-1/2}$ for the values of $Z$ of interest,
including the set $|Z|=1$ for which the sequential action of intertwining
operators is defined.
Under this restriction, we can replace $\Phi(x\pm\eta)$ in the above equation
with $\Phi(x)$ after appropriate shifts of the arguments in integrands
and reduce the integrations to integrations over the interval $[0,1]$.
Because $\Phi(x)$ is
an arbitrary function, we obtain the equation
\begin{eqnarray*}
&& \frac{\Gamma(\pm z\pm (x-\eta) -g,\pm 2x)}
{\Gamma(\pm z\pm x-g,\pm 2(x-\eta))}
\frac{e^{2\pi \textup{i} (\eta-x)} }{ \mathrm{R}(\tau)\theta_1(2x-2\eta)}\bar \theta_k(x-g-\eta)
\\ && \makebox[2em]{}
-\frac{\Gamma(\pm z\pm (x+\eta) -g,\pm 2x)}
{\Gamma(\pm z\pm x-g,\pm 2(x+\eta))}
\frac{e^{2\pi \textup{i} (\eta+x)}}{\mathrm{R}(\tau)\theta_1(2x+2\eta)}\bar \theta_k(x+g+\eta)
=\frac{\Gamma(-2g)}{\Gamma(-2g-2\eta)}\bar \theta_k(z)
\end{eqnarray*}
for the integrands.
Appropriately simplifying the elliptic gamma function ratios,
we obtain precisely identity \eqref{id2}. Hence, second operator
contiguous relation (see \eqref{RR2}) holds under milder
conditions than are necessary for inversion relation \eqref{inv} to hold.

\section{The intertwining operator for two-index discrete lattices
$g = n\eta+m\frac{\tau}{2}$ and $g = \frac{1}{2} +n\eta+m\frac{\tau}{2}$
with $n,m \in \mathbb{Z}_{\geq 0}$}

Contiguous relations for the intertwining operator lead to particular factored
forms of  $\mathrm{M}(g)$ at special integer-valued
points  of the spin lattice $g = n\eta+m\frac{\tau}{2}$
and $g =  \frac{1}{2}+n\eta+m\frac{\tau}{2}$, where $n,m\in\Z_{\geq0}$.

We can repeat the above considerations with the change
$2\eta \rightleftarrows\tau$ and obtain two types of the recurrence relations  ($k=3,4$)
\begin{equation}
\mathrm{A}_k(g)\,\mathrm{M}(g) = \mathrm{M}(g+\eta)\,\theta_k
 \left(z | {\textstyle\frac{\tau}{2}}\right)\, , \ \ \
\mathrm{B}_k(g)\,\mathrm{M}(g) = \mathrm{M}\left(g+{\textstyle\frac{\tau}{2}}\right)\,
\theta_k \left(z | \eta\right)\,,
\label{RR}\end{equation}
where $\mathrm{A}_k(g)$ and $\mathrm{B}_k(g)$ are the difference operators
$$
\mathrm{A}_k(g) = e^{\pi \textup{i}\frac{z^2}{ \eta}}\,\frac{c_A}{\theta_1(2z | \tau)}
\left[ \theta_k \left(z+g+\eta | {\textstyle\frac{\tau}{2}}\right)\, e^{\eta \partial_z} -
\theta_k \left(z-g-\eta | {\textstyle\frac{\tau}{2}}\right)\, e^{-\eta \partial_z}
\right]\, e^{-\pi \textup{i}\frac{z^2}{ \eta}}\,,
$$
$$
\mathrm{B}_k(g) = e^{2\pi \textup{i}\frac{z^2}{ \tau}}
\,\frac{c_B}{\theta_1(2z | 2\eta)}
\left[ \theta_k \left(z+g+{\textstyle\frac{\tau}{2}}| \eta\right)
\, e^{{\textstyle\frac{\tau}{2}} \partial_z} -
\theta_k \left(z-g-{\textstyle\frac{\tau}{2}}| \eta\right)
\, e^{-{\textstyle\frac{\tau}{2}} \partial_z}
\right]\, e^{-2\pi \textup{i}\frac{z^2}{ \tau}}
$$
and
$$
c_A = \frac{e^{\pi \textup{i} \eta}}{\mathrm{R}(\tau)}\ \ \,,\ \
c_B = \frac{e^{\pi \textup{i} {\textstyle\frac{\tau}{2}}}}{\mathrm{R}(2\eta)}\,.
$$

Using the initial condition $\mathrm{M}(0) = \II$, which is proved by simple
residue calculus \cite{DS}, we can solve the recurrence relations
and obtain
$$
\mathrm{M}^{(k)}(n\eta) = \mathrm{A}_k(n\eta-\eta)\cdots \mathrm{A}_k(\eta) \mathrm{A}_k(0)\cdot
\theta_k^{-n} \left(z | {\textstyle\frac{\tau}{2}}\right)\,,
$$
$$
\mathrm{M}^{(k)}\left( m {\textstyle\frac{\tau}{2}}\right) = \mathrm{B}_k\left(m {\textstyle\frac{\tau}{2}} -{\textstyle\frac{\tau}{2}}\right)\cdots \mathrm{B}_k\left({\textstyle\frac{\tau}{2}}\right) \mathrm{B}_k(0)\cdot
\theta_k^{-m} \left(z | \eta\right)\,.
$$
Although the form of the intertwining operator should be independent of $k$,
$\mathrm{M}^{(k)}(g)\equiv \mathrm{M}(g)$,
we introduced an additional upper index $k$ in order to indicate
a possible dependence on it.

The derived expressions are particular cases of the general operator
$\mathrm{M}^{(k)}\left(n\eta + m {\textstyle\frac{\tau}{2}}\right)$.
But they are used as building blocks for constructing this general operator.
We first transform these factored operators to finite sums with explicit coefficients.
As an illustration, we consider two simple examples:
\begin{eqnarray}\nonumber &&
e^{-\pi \textup{i}\frac{z^2}{ \eta}}\,\mathrm{M}^{(k)}(\eta)\,e^{\pi \textup{i}\frac{z^2}{ \eta}}\, = \frac{c_A}{\theta_1(2z)}
\left[\theta_k(z+\eta|{\textstyle\frac{\tau}{2}}) \,e^{\eta \partial_z} -
\theta_k(z-\eta|{\textstyle\frac{\tau}{2}}) \, e^{-\eta \partial_z}\right]\,
\theta_k^{-1} \left(z | {\textstyle\frac{\tau}{2}}\right) =
\\ && \makebox[6em]{}
= \frac{c_A}{\theta_1(2z)}
\left[e^{\eta \partial_z} - e^{-\eta \partial_z}\right] \,,
\label{meta}\end{eqnarray}
\begin{eqnarray*} &&
e^{-\pi \textup{i}\frac{z^2}{ \eta}}\,\mathrm{M}^{(k)}(2 \eta)\,e^{\pi \textup{i}\frac{z^2}{ \eta}} = \frac{c_A^2}{\theta_1(2z)}
\left[\theta_k(z+2\eta|{\textstyle\frac{\tau}{2}}) \,e^{\eta \partial_z} -
\theta_k(z-2\eta|{\textstyle\frac{\tau}{2}}) \, e^{-\eta \partial_z}\right]\cdot
\\ && \makebox[6em]{}
\cdot\,\frac{1}{\theta_1(2z)}
\left[e^{\eta \partial_z} - e^{-\eta \partial_z}\right]\,
\theta_k^{-1} \left(z | {\textstyle\frac{\tau}{2}}\right)\,.
\end{eqnarray*}
Expanding the latter expression we obtain a sum of four finite-difference
operators, which can be transformed to the form
$$
e^{-\pi \textup{i}\frac{z^2}{ \eta}}\,\mathrm{M}^{(k)}(2 \eta)\,e^{\pi \textup{i}\frac{z^2}{ \eta}} =
\frac{c_A^2}{\theta_1(2 z-2\eta)\theta_1(2 z)\theta_1(2 z+2\eta)}\, \cdot
$$
$$
\cdot\left[
\theta_1(2 z-2\eta)\, e^{2 \eta \partial_z}
- \frac{\theta_1(4 \eta)}{\theta_1(2 \eta)} \theta_1(2 z)
+\theta_1(2 z+2\eta)\,e^{-2 \eta \partial_z}
\right]
$$
using identity (\ref{theta}).
In these examples, we can drop the index $k$ in the
notation $\mathrm{M}^{(k)}(n \eta)$ because the result is independent
of $k$, as expected.

Take the general ansatz for the intertwining operator,
$$
e^{-\pi \textup{i}\frac{z^2}{ \eta}}\,\mathrm{M}(n \eta)\,
e^{\pi \textup{i}\frac{z^2}{ \eta}} = \sum_{\ell=0}^{n}
(-1)^\ell\, \alpha_\ell^{(n)}(z)\,
e^{(n-2\ell)\eta \partial_z}\, ,
$$
and substitute it in the equality  $\mathrm{M}((n+1) \eta)
=\mathrm{A}_k(n\eta)\mathrm{M}(n \eta)\theta_k(z|{\textstyle\frac{\tau}{2}})^{-1}$.
For $\ell=1,\ldots,n$ this yields the recurrence relations
\begin{equation}
\alpha_\ell^{(n+1)}(z)=\frac{c_A}{\theta_1(2z|\tau)}\cdot
\frac{\theta_k(z+(n+1)\eta|{\textstyle\frac{\tau}{2}})\alpha_\ell^{(n)}(z+\eta)
+\theta_k(z-(n+1)\eta|{\textstyle\frac{\tau}{2}})\alpha_{\ell-1}^{(n)}(z-\eta)}
{\theta_k(z+(n+1-2\ell)\eta|{\textstyle\frac{\tau}{2}}) }.
\label{recrel}\end{equation}
For $\ell=0$ and $\ell=n+1$, we have
$$
\alpha_0^{(n+1)}(z)=\frac{c_A}{\theta_1(2z|\tau)}\alpha_0^{(n)}(z+\eta), \quad
\alpha_{n+1}^{(n+1)}(z)=\frac{c_A}{\theta_1(2z|\tau)}\alpha_n^{(n)}(z-\eta).
$$
Because $\alpha_0^{(0)}(z)=1$, we obtain
$$
\alpha_0^{(n)}(z)=\frac{c_A^n}{\prod_{k=0}^{n-1}\theta_1(2z+2\eta k|\tau)},
\quad
\alpha_n^{(n)}(z)=\frac{c_A^n}{\prod_{k=0}^{n-1}\theta_1(2z-2\eta k|\tau)}.
$$
These boundary values for $\alpha_{\ell}^{(n)}(z)$ uniquely fix the solution of recurrence
relation \eqref{recrel}, which can be found by induction:
\begin{eqnarray}\label{alpha} &&
\alpha_\ell^{(n)}(z) = c_A^n \cdot\begin{bmatrix}
n \\ \ell
\end{bmatrix}_{\tau,2\eta}\cdot
\frac{\theta_1(2 z + 2 \eta(n-2\ell)\,|\tau)}{\prod_{j = 0}^{n}
\theta_1(2 z - 2\eta(\ell-j)\, |\tau)}\ \ \,,\ \ \
\\ &&  \makebox[2em]{}
\begin{bmatrix}
n \\ \ell
\end{bmatrix}_{\tau,2\eta} =
\frac{\prod_{j=1}^{n}
\theta_1(2 \eta j\,|\tau)}
{\prod_{j=1}^{\ell} \theta_1(2 \eta j\,|\tau)
\cdot \prod_{j=1}^{n-\ell} \theta_1(2 \eta j\,|\tau)}\,.
\nonumber\end{eqnarray}
In the expression for $\alpha_\ell^{(n)}(z)$, we extract
the elliptic binomial coefficient $\begin{bmatrix}
n \\ \ell
\end{bmatrix}_{\tau,2\eta}$, which is independent of $z$, and
explicitly write the remaining $z$-depended part.
The discrete intertwining operator $\mathrm{M}(n\eta)$
was first obtained in \cite{Z} in this form. As we see,
indeed, the general result has no $k$-index dependence, which is
present in the recurrence relation for $\alpha_\ell^{(n)}(z)$.

A similar expression for the operator
$\mathrm{M}\left( m {\textstyle\frac{\tau}{2}}\right)$ is
obtained by a simple permutation $c_A \rightleftarrows c_B$,
$n \rightleftarrows m$, and $\tau \rightleftarrows 2\eta$:
$$
e^{-2\pi \textup{i}\frac{z^2}{ \tau}}\,\mathrm{M}\left( m {\textstyle\frac{\tau}{2}}\right)\,e^{2\pi \textup{i}\frac{z^2}{\tau}} =
\sum_{\ell=0}^{m}
(-1)^\ell\, \mathrm{\beta}^{(m)}_{\ell}(z)\,
e^{(m-2\ell)\,\frac{\tau}{2} \partial_z}\,,
$$
where
\begin{eqnarray}\label{beta} &&
\mathrm{\beta}^{(m)}_{\ell}(z) =
c_B^m \cdot
\begin{bmatrix}
m \\ \ell
\end{bmatrix}_{2\eta,\tau}\cdot
\frac{\theta_1(2 z + \tau(m - 2\ell)\,|2\eta)}
{\prod_{j= 0}^{m}\theta_1(2 z - \tau(\ell -j)\, |2\eta)}
\ \ \,,
\\ &&  \makebox[2em]{}
\ \ \begin{bmatrix} m \\ \ell \end{bmatrix}_{2\eta,\tau} =
\frac{\prod_{j=1}^{m}
\theta_1(\tau j\,|2\eta)}
{\prod_{j=1}^{\ell} \theta_1(\tau j\,|2\eta)
\cdot \prod_{j=1}^{m-\ell} \theta_1(\tau j\,|2\eta)}\,.
\nonumber\end{eqnarray}

We now describe the general case. It is easy to derive the representation
\begin{eqnarray}\nonumber  &&
\mathrm{M}^{(k)}\left(n\eta + m {\textstyle\frac{\tau}{2}}\right) =
\mathrm{A}_k(n\eta-\eta+ m {\textstyle\frac{\tau}{2}})\cdots
\mathrm{A}_k(\eta+m {\textstyle\frac{\tau}{2}})
\mathrm{A}_k(m {\textstyle\frac{\tau}{2}})\,\cdot
\\  && \makebox[2em]{}
\cdot\,
\mathrm{B}_k\left(m {\textstyle\frac{\tau}{2}} -{\textstyle\frac{\tau}{2}}\right)\cdots \mathrm{B}_k\left({\textstyle\frac{\tau}{2}}\right) \mathrm{B}_k(0)
\cdot
\theta_k^{-m} \left(z | \eta\right)
\theta_k^{-n} \left(z | {\textstyle\frac{\tau}{2}}\right)\,.
\label{genform}\end{eqnarray}
Of course, there are many equivalent ways to represent
$\mathrm{M}^{(k)}\left(n\eta + m {\textstyle\frac{\tau}{2}}\right)$ as a product of
$\mathrm{A}_k$ and $\mathrm{B}_k$ operators, and we have described only one of them.

We can consider the lattice $g={\textstyle\frac{1}{2}}+n\eta +
m {\textstyle\frac{\tau}{2}}$ analogously. We must use the fact that
$\mathrm{M}({\textstyle\frac{1}{2}})=P$,
where $P=e^{\frac{1}{2}\partial_z}$ is the operator of shifting $z$ by 1/2 \cite{DS}.
Because we work in the space of functions that are meromorphic in
$w=e^{2\pi\textup{i}z}$, the half period shifting is equivalent to the parity
transformation, $Pw=-w$. Repeating the previous procedure once more in this setting,
we obtain
\begin{eqnarray}\nonumber &&
\mathrm{M}^{(k)}\left(n\eta + m {\textstyle\frac{\tau}{2}}+{\textstyle\frac{1}{2}}\right) =
\mathrm{A}_k(n\eta-\eta+ m {\textstyle\frac{\tau}{2}}+{\textstyle\frac{1}{2}})\cdots
\mathrm{A}_k(\eta+m {\textstyle\frac{\tau}{2}}+{\textstyle\frac{1}{2}})
\mathrm{A}_k(m {\textstyle\frac{\tau}{2}}+{\textstyle\frac{1}{2}})\,\cdot
\\  && \makebox[2em]{}
\cdot\,
\mathrm{B}_k\left(m {\textstyle\frac{\tau}{2}} -{\textstyle\frac{\tau}{2}}
+{\textstyle\frac{1}{2}}\right)\cdots \mathrm{B}_k\left({\textstyle\frac{\tau}{2}}
+{\textstyle\frac{1}{2}}\right) \mathrm{B}_k({\textstyle\frac{1}{2}})
\cdot
\theta_k^{-m} \left(z+ \textstyle{\frac{1}{2}} | \eta\right)
\theta_k^{-n} \left(z+ \textstyle{\frac{1}{2}} | {\textstyle\frac{\tau}{2}}\right)P\,.
\label{S1/2}\end{eqnarray}

Because $\theta_3(z\pm{\textstyle\frac{1}{2}})=\theta_4(z)$, we see that
$$
\mathrm{A}_{3,4}(g+{\textstyle\frac{1}{2}})=\mathrm{A}_{4,3}(g), \qquad
\mathrm{B}_{3,4}(g+{\textstyle\frac{1}{2}})=\mathrm{B}_{4,3}(g).
$$
Therefore, the intertwining operator for the second lattice is obtained
from the first simply by exchanging $\theta_3$ and $\theta_4$ and
multiplying by $P$ from the right:
\begin{equation}
\mathrm{M}^{(3,4)}\left(n\eta + m {\textstyle\frac{\tau}{2}}
+{\textstyle\frac{1}{2}}\right)
=\mathrm{M}^{(4,3)}\left(n\eta + m {\textstyle\frac{\tau}{2}}\right)P.
\label{inter1/2}\end{equation}

\section{A finite-dimensional invariant space of the elliptic modular double}

The irreducible representation of the Sklyanin algebra
for $g=n\eta$ at (half)-integer spin $\ell=\frac{n-1}{2}$, $n\in\Z_{> 0}$,
is $n$-dimensional and can be realized in the space $\Theta^{+}_{2n-2}(z|\tau)$
consisting of even theta functions of $z$ of order $2n-2$
with the quasiperiods 1 and $\tau$ \cite{skl2}. Let
$$
\mathrm{W}(g)=e^{-\pi\textup{i}z^2/\eta}\,\mathrm{M}(g)\,
e^{\pi\textup{i}z^2/\eta}.
$$
From the factored representation of the operator $\mathrm{W}(\eta n)$,
we see that its action annihilates the $\Theta^{+}_{2n-2}(z|\tau)$ space elements.
We demonstrate this fact using the recurrence relations.

The function $\theta_k^2(z|\tau)$ for any $k=1,\ldots,4$
is an even theta function of the second order. We fix two indices
$k_1$ and $k_2,\, k_1\neq k_2,$ and set $e_1:=\theta_{k_1}^2(z|\tau)$ and
$e_2=\theta_{k_2}^2(z|\tau)$.
All monomials of these building blocks $e_1^k\,e^{N-k}_2$,
$k=0,\ldots,N$, are even theta functions of order $2N$. They are linearly
independent, i.e., the equality $\sum_{k=0}^Nc_ke_1^k(z)\,e^{N-k}_2(z)=0$
is satisfied only for $c_k=0$ (to verify this, we must
sequentially set $z$ equal to a root of the equation
$\theta_{k_1}(z|\tau)=0$ or $\theta_{k_2}(z|\tau)=0$ after taking an appropriate
number of derivatives with respect to $z$).
The number of such monomials yields the dimension of the space of even
theta functions of order $2N$, and these monomials therefore form a basis
of $\Theta_{2N}^+(z|\tau)$.

Using the formulas
\begin{eqnarray}\nonumber &&
2\,\theta_1(x+y)\,\theta_1(x-y) = \bar\theta_4(x)\,\bar\theta_3(y)
-\bar\theta_4(y)\,\bar\theta_3(x),
\\ \nonumber &&
2\,\theta_2(x+y)\,\theta_2(x-y) = \bar\theta_3(x)\,\bar\theta_3(y)
-\bar\theta_4(y)\,\bar\theta_4(x),
\\ \label{33} &&
  2\,\theta_3(x+y)\,\theta_3(x-y)
= \bar\theta_3(x)\,\bar\theta_3(y) +\bar\theta_4(y)\,\bar\theta_4(x),
\\ \nonumber &&
2\, \theta_4(x+y)\,\theta_4(x-y)=\bar\theta_4(x)\,\bar\theta_3(y)
+\bar\theta_4(y)\,\bar\theta_3(x),
\end{eqnarray}
we can express any square $\theta_k^2(z|\tau)$ as a
linear combination of $\bar\theta_3(z)$ and $\bar\theta_4(z)$.
Therefore, the set of $n$ functions
$$
\left[\theta_3 \left(z | {\textstyle\frac{\tau}{2}}\right)\right]^j \,
\left[\theta_4 \left(z | {\textstyle\frac{\tau}{2}}\right)\right]^{n-1-j}
\;\; , \;\; j=0,1,\ldots, n-1,
$$
forms a basis in the space $\Theta^{+}_{2n-2}(z|\tau)$. We note that
these functions can be multiplied by an arbitrary even theta function
$e^{\alpha z^2}$, $\alpha\in\mathbb{C}$, of zeroth order because the results are
even theta functions as before with some quasiperiodicity exponential multipliers
appearing after the shifts of $z$ by $1$ and $\tau$. The value of the parameter
$\alpha$ fixing a concrete form of the invariant space of theta functions
depends on the realization of the Sklyanin algebra generators. The standard form
of generators suggested in \cite{skl2} corresponds to the choice $\alpha=0$.

Applying the recurrence relation $n-1$ times and taking the relation
\begin{equation}
\mathrm{W}(\eta)\cdot 1 = \frac{c_A}{\theta_1(2z)}
\left[e^{\eta \partial_z} - e^{-\eta \partial_z}\right]
\cdot 1 = 0
\label{quasic}\end{equation}
into account (here 1 can be replaced by an arbitrary $2\eta$-periodic function),
for $n>0$, we obtain
$$
\mathrm{W}(\eta n) \cdot \left[\theta_3 \left(z | {\textstyle\frac{\tau}{2}}\right)\right]^j \,
\left[\theta_4 \left(z | {\textstyle\frac{\tau}{2}}\right)\right]^{n-1-j}  =0\,.
$$
Hence, zero modes of the operator  $\mathrm{M}(n\eta)$ for $n>0$ have the form
$$
\mathrm{M}(n\eta ) \cdot \phi_j(z)\psi_j(z)=0,\quad
\phi_j(z):=
e^{\pi\textup{i}z^2/\eta}\left[\theta_3 \left(z | {\textstyle\frac{\tau}{2}}\right)\right]^j \,
\left[\theta_4 \left(z | {\textstyle\frac{\tau}{2}}\right)\right]^{n-1-j},
$$
where $j=0,1,\ldots, n-1$ and $\psi_j(z)$ is an arbitrary periodic function,
$\psi_j(z+2\eta)=\psi_j(z)$.

The usual basis of a finite-dimensional representation of the Sklyanin algebra
at (half-)integer values of the spin $\ell$ does not contain the exponential
multiplier $e^{\pi\textup{i}z^2/\eta}$. This multiplier emerged because of the
non-standard  form of the Sklyanin algebra generators \eqref{Sklyan}
we use, and it does not change the space $\Theta^{+}_{2n-2}$,
being an even theta function of zeroth order. But we note that it is
not analytic in the variable $e^{2\pi \textup{i} z}$, which contrasts
with the analytical structure of the intertwining operator $\mathrm{M}(g)$.

We consider the invariance of the described space of zero modes under the
action of generators of the elliptic modular double. For $\mathbf{S}^a(n\eta)$,
we have
\begin{eqnarray} \nonumber && \makebox[-1em]{}
\mathbf{S}^a(n\eta) \phi_j(z)\psi_j(z) =
\psi_j(z+\eta) e^{\pi\textup{i}\frac{z^2}{\eta}}
\frac{\textup{i}^{\delta_{a,2}}\theta_{a+1}(\eta|\tau)}{\theta_1(2 z|\tau) }
\cdot \\  \nonumber && \makebox[0em]{} \cdot
\Bigl[\,\theta_{a+1} \left(2z-n\eta+\eta|\tau\right)e^{\eta\partial_z}
- \theta_{a+1}
\left(-2z-n\eta+\eta|\tau\right)e^{-\eta\partial_z}\Bigl]
\left[\theta_3 \left(z | {\textstyle\frac{\tau}{2}}\right)\right]^j \,
\left[\theta_4 \left(z | {\textstyle\frac{\tau}{2}}\right)\right]^{n-1-j}.
\end{eqnarray}
It is easy to see that the action of the finite-difference operator in the right-hand side
leads to theta functions of the order $2n-2$ (for this, it suffices to
verify the holomorphicity of this function and its double quasiperiodicity).
Therefore, for a fixed function $\psi_j(z)=\psi(x)$ satisfying the
half (anti)periodicity condition $\psi(z+\eta)=\pm \psi(z)$, we obtain
a finite-dimensional invariant space for the usual Sklyanin algebra
found in \cite{skl2}. But this space is not invariant under the action
of generators of the modular pair $\mathbf{\tilde S}^a$. Indeed, we have
$$
\mathbf{\tilde S}^a(n\eta)\phi_j(z)\psi(z)=
(-1)^{n-1-j}e^{\pi\textup{i}(\frac{\tau^2}{4\eta} -n\frac{\tau}{2})}\phi_j(z)
\textup{i}^{\delta_{a,2}}\theta_{a+1}({\textstyle\frac{\tau}{2}}|2\eta)\,
D^a\psi(z),
$$
where
\begin{equation} \label{Da}
D^a=f_a(z)e^{{\textstyle\frac{\tau}{2}}\partial_z}
+f_a(-z)  e^{-{\textstyle\frac{\tau}{2}}\partial_z}
\end{equation}
and
$$
f_a(z):=\frac{e^{\pi\textup{i}(\frac{\tau}{\eta}-2n)z}\theta_{a+1}
\left(2z-n\eta+{\textstyle\frac{\tau}{2}}|2\eta\right)}{\theta_1(2 z|2\eta)}.
$$
Hence, we obtain an invariant space if $\psi(z)$ is
an eigenfunction of all four operators $D^a$ simultaneously,
$D^a\psi(z)=\lambda_a\psi(z)$, which is impossible. Indeed, we have
$f_{0,3}(z+\eta)=f_{0,3}(z)$ and $f_{1,2}(z+\eta)=-f_{1,2}(z)$.
Because $\psi(z+\eta)=\pm \psi(z)$, we conclude that $\lambda_{1,2}=0$,
whence we obtain a contradiction:
$$
\frac{\psi(z+{\textstyle\frac{\tau}{2}}) }{ \psi(z-{\textstyle\frac{\tau}{2}}) }
=-\frac{f_2(-z)}{f_2(z)} = -\frac{f_1(-z)}{f_1(z)}.
$$

We now choose the $2\eta$-periodic functions mentioned above in the form
$$
\psi_j(z)=\frac{e^{-\pi\textup{i}z^2/\eta}\psi(z)}{\theta(q^{-1/4}Z^{\pm 1};q)},
\quad \theta(tZ^{\pm1};q):= \theta(tZ;q)\theta(tZ^{-1};q),
$$
and obtain
\begin{equation}
\mathrm{M}(n\eta ) \varphi_{j}^{(n)}(z)\psi(z)=0,\quad
 \varphi_{j}^{(n)}(z):=
\frac{\left[\theta_3 \left(z | {\textstyle\frac{\tau}{2}}\right)\right]^j \,
\left[\theta_4 \left(z | {\textstyle\frac{\tau}{2}}\right)\right]^{n-1-j}}
{\theta(q^{-1/4}Z^{\pm 1};q)}, \quad  j=0,\ldots, n-1.
\label{zmneta}\end{equation}
In the first case, zero modes are holomorphic in $z$
but they are not analytic functions of $Z$. Functions
\eqref{zmneta} are meromorphic in $Z$, i.e., they have a qualitatively
different analytic nature. Requiring analyticity in $Z$
leads to reducing the gauge freedom: for an arbitrary
multiplier $\psi(z)=\psi(z+2\eta)$, we now have an additional constraint
$\psi(z+1)=\psi(z)$. This means that $\varphi_{j}^{(n)}(z)$
can be multiplied by an elliptic function $\chi(Z)$ such that $\chi(qZ)=\chi(Z)$.

We consider the invariance of the described space of zero modes. We have
\begin{eqnarray} \nonumber && \makebox[-1em]{}
\mathbf{S}^a(n\eta) \varphi_{j}^{(n)}(z)\chi(Z)=
-\frac{\chi(q^{1/2}Z)}{\theta(q^{-1/4}Z^{\pm 1};q)}
\frac{\textup{i}^{\delta_{a,2}}\theta_{a+1}(\eta|\tau)}{\theta_1(2 z|\tau) }
\cdot \\  && \makebox[0em]{}
\cdot \Bigl[\,\theta_{a+1} \left(2z-n\eta+\eta|\tau\right)e^{\eta\partial_z}
- \theta_{a+1}
\left(-2z-n\eta+\eta|\tau\right)e^{-\eta\partial_z}\Bigl]
\left[\theta_3 \left(z | {\textstyle\frac{\tau}{2}}\right)\right]^j \,
\left[\theta_4 \left(z | {\textstyle\frac{\tau}{2}}\right)\right]^{n-1-j}.
\nonumber\end{eqnarray}
The action of the finite-difference operator in the right-hand side remains
in the space $\Theta^+_{2n-2}(z|\tau)$. We therefore have an
invariant space if $\chi(q^{1/2}Z)=\pm \chi(Z)$.

The action of the operators  $\mathbf{\tilde S}^a(n\eta)$ has the form
$$
\mathbf{\tilde S}^a(n\eta)\varphi_{j}^{(n)}(z)\chi(Z)=
(-1)^{n-1-j}e^{-\pi\textup{i}n\frac{\tau}{2}}
\varphi_{j}^{(n)}(z)\,
\textup{i}^{\delta_{a,2}}\theta_{a+1}({\textstyle\frac{\tau}{2}}|2\eta)\,D^a\chi(Z)
$$
where $D^a$ is given by formula \eqref{Da} and in this case
$$
f_a(z):=e^{-2\pi\textup{i}nz}\frac{\theta_{a+1}
\left(2z-n\eta+{\textstyle\frac{\tau}{2}}|2\eta\right)\theta(q^{-1/4}Z^{\pm 1};q)}
{\theta_1(2 z|2\eta)\theta(q^{-1/4}(p^{1/2}Z)^{\pm1};q)}.
$$
As above, $f_{0,3}(z+\eta)=f_{0,3}(z)$ and $f_{1,2}(z+\eta)=-f_{1,2}(z)$.
Correspondingly, there is no function  $\chi(Z)$ that
would be a simultaneous eigenfunction of four operators $D^a$.

Therefore, for the spin values $g=n\eta, {\textstyle\frac{1}{2}}+n\eta,
m{\textstyle\frac{\tau}{2}},$ and ${\textstyle\frac{1}{2}}+m{\textstyle\frac{\tau}{2}}$,
we obtain finite-dimensional invariant spaces only for one of the Sklyanin algebras
in the elliptic modular double.

We now consider other discrete values of the spin $g$. Using relations
\eqref{RR}, we can write
$$
\mathrm{M}(n\eta+\textstyle{\frac{\tau}{2}})=\mathrm{B}_k(n\eta)\mathrm{M}(n\eta )
\theta_k^{-1}(z|\eta).
$$
We move the factor $\theta_k^{-1}(z|\eta)$ to the left
using its quasiperiodicity and the explicit expression for
$\mathrm{M}(n\eta )$ and obtain
$$
\mathrm{M}(n\eta+\textstyle{\frac{\tau}{2}})=\mathrm{\tilde B}(n\eta)\, \mathrm{W}(n\eta ),
$$
where
\begin{eqnarray*} &&
\mathrm{\tilde B}(n\eta)= e^{2\pi \textup{i}\frac{z^2}{ \tau}}
\,\frac{c_B(-1)^{n\delta_{k,4}}}{\theta_1(2z | 2\eta)}
\left[ \frac{\theta_k \left(z+n\eta+{\textstyle\frac{\tau}{2}}| \eta\right)}
{\theta_k \left(z+{\textstyle\frac{\tau}{2}}| \eta\right)}
\, e^{{\textstyle\frac{\tau}{2}} \partial_z} -
\frac{\theta_k \left(z-n\eta-{\textstyle\frac{\tau}{2}}| \eta\right)}
{\theta_k \left(z-{\textstyle\frac{\tau}{2}}| \eta\right)}
\, e^{-{\textstyle\frac{\tau}{2}} \partial_z}
\right]\, e^{-2\pi \textup{i}\frac{z^2}{ \tau}}\,.
\end{eqnarray*}
It is easy to see that the ratio of theta functions in the operator
$\mathrm{\tilde B}$ can be simplified and that $\mathrm{\tilde B}$
is proportional to $e^{\tau\partial_z/2}-e^{-\tau\partial_z/2}$
conjugated with exponential multipliers. Hence,
$$
\mathrm{M}(n\eta+\textstyle{\frac{\tau}{2}})\cdot \left[\theta_3 \left(z | {\textstyle\frac{\tau}{2}}\right)\right]^j \,
\left[\theta_4 \left(z | {\textstyle\frac{\tau}{2}}\right)\right]^{n-1-j}  =0\, ,
\; j=0,1,\ldots, n-1,\ n>0,
$$
i.e., the same functions as before but without the additional
(exponential or theta functional) factor form zero modes of the more
complicated operator $\mathrm{M}(n\eta+\textstyle{\frac{\tau}{2}})$. For $n=0$,
we have $\mathrm{W}(0)=\II$, and the zero mode of
$\mathrm{M}(\textstyle{\frac{\tau}{2}})$ is therefore determined by the factor
$\mathrm{\tilde B}(0)$, and the requirement of
analyticity in $Z$ obviously leads to the function
$\chi(Z)/\theta(p^{-1/4}Z^{\pm1};p)$, $\chi(pZ)=\chi(Z)$.

In the same way, permuting the parameters $\tau$ and $2\eta$, we find
$$
\mathrm{M}\left(m{\textstyle\frac{\tau}{2}} \right)
\cdot e^{2\pi\textup{i}z^2/\tau}\left[\theta_3 \left(z | \eta \right)\right]^\ell \,
\left[\theta_4 \left(z | \eta \right)\right]^{m-1-\ell}\psi(z)  =0\, ,
 \  m>0\, , \quad \psi(z+{\textstyle\frac{\tau}{2}})=\pm\psi(z),
$$
where $\ell=0,1,\cdots, m-1$, $m>0$, or
$$
\mathrm{M}\left(m{\textstyle\frac{\tau}{2}} \right) \cdot
\frac{\left[\theta_3 \left(z |\eta\right)\right]^\ell \,
\left[\theta_4 \left(z | \eta\right)\right]^{m-1-\ell}}
{\theta(p^{-1/4}Z^{\pm 1};p)}\chi(Z)=0\, , \quad  \chi(p^{1/2}Z)=\pm\chi(Z).
$$

For $\mathrm{M}\left(m{\textstyle\frac{\tau}{2}}+\eta\right)$, we obtain
$$
\mathrm{M}\left(m{\textstyle\frac{\tau}{2}}+\eta\right)
\cdot \left[\theta_3 \left(z | \eta \right)\right]^\ell \,
\left[\theta_4 \left(z | \eta \right)\right]^{m-1-\ell}  =0\, , \
\ell=0,1,\ldots, m-1, \  m>0.
$$

As is shown below, the set of $nm$ functions
\begin{equation}
\varphi_{j,\ell}^{(n,m)}(z):=
\left[\theta_3 \left(z | {\textstyle\frac{\tau}{2}}\right)\right]^j \,
\left[\theta_4 \left(z | {\textstyle\frac{\tau}{2}}\right)\right]^{n-1-j}
\cdot
\left[\theta_3 \left(z | \eta \right)\right]^\ell \,
\left[\theta_4 \left(z | \eta \right)\right]^{m-1-\ell}, \, \  n,m>0,
\label{varphi}\end{equation}
is annihilated by the operator
$\mathrm{M}\left(n\eta+m{\textstyle\frac{\tau}{2}}\right)$,
\begin{equation}
\mathrm{M} \left(n\eta+m{\textstyle\frac{\tau}{2}}\right)\cdot
\varphi_{j,\ell}^{(n,m)}(z)= 0\, ,  \quad  n,m>0.
\label{Mphi}\end{equation}
We can multiply $\varphi_{j,\ell}^{(n,m)}(z)$ by arbitrary periodic functions
with either the period $2\eta$ or $\tau$, and they are still zero modes
of $\mathrm{M}\left(n\eta+m{\textstyle\frac{\tau}{2}}\right)$.

The operator $\mathrm{M}({\textstyle\frac{1}{2}}+n\eta + m {\textstyle\frac{\tau}{2}})$
has analogous properties because the shift of $z$ by 1/2
simply permutes $\theta_3(z)$ and $\theta_4(z)$, i.e., we have a basis shuffle,
$$
\mathrm{M} \left({\textstyle\frac{1}{2}}+n\eta+m{\textstyle\frac{\tau}{2}}\right)
\cdot \varphi_{n-1-j,m-1-\ell}^{(n,m)}(z)= 0\, , \quad n,m>0.
$$

The intertwining operator kernel forms an invariant space under the action of
the Sklyanin algebra generators. But we treat only a finite-dimensional
subspace of this kernel. We must therefore verify that it is
invariant. For this, we apply the operators $\mathbf{S}^a(n\eta+m\tau/2)$
to $\varphi_{j,\ell}^{(n,m)}(z)$. Because the functions
$\theta_{3,4} \left(z | \eta \right)$ are quasiperiodic under the shifts of  $z$
by $\eta$, we have
\begin{eqnarray} \nonumber && \makebox[-2em]{}
\mathbf{S}^a(n\eta+m\tau/2)\varphi_{j,\ell}^{(n,m)}(z)=
(-1)^{m-1-\ell}e^{-\pi\textup{i}m\eta}
\left[\theta_3 \left(z | \eta \right)\right]^\ell \,
\left[\theta_4 \left(z | \eta \right)\right]^{m-1-\ell}
\frac{\textup{i}^{\delta_{a,2}}
\theta_{a+1}(\eta|\tau)}{\theta_1(2 z|\tau) } \cdot
\\  \nonumber &&
\cdot \Bigl[\,e^{-2\pi\textup{i}mz} \theta_{a+1} \left(2
z-n\eta-m\tau/2 +\eta|\tau\right)e^{\eta\partial_z} -
\\  &&
- e^{2\pi\textup{i}mz}\theta_{a+1}
\left(-2z-n\eta-m\tau/2+\eta|\tau\right)e^{-\eta\partial_z}\Bigl]
\left[\theta_3 \left(z | {\textstyle\frac{\tau}{2}}\right)\right]^j \,
\left[\theta_4 \left(z | {\textstyle\frac{\tau}{2}}\right)\right]^{n-1-j}.
\nonumber\end{eqnarray}

Let the integer $m$ be even. In the arguments of the theta functions
$\theta_{a+1}$, there are then shifts by integer multiples of $\tau$.
Using the qusiperiodicity of $\theta_{a+1}$, we obtain
\begin{eqnarray} \nonumber && \makebox[-2em]{}
\mathbf{S}^a(n\eta+m\tau/2)\varphi_{j,\ell}^{(n,m)}(z)=\mu_a\,
e^{-\pi\textup{i}mn\eta-\pi\textup{i}m^2\tau/4}
\left[\theta_3 \left(z | \eta \right)\right]^\ell \,
\left[\theta_4 \left(z | \eta \right)\right]^{m-1-\ell}
\frac{\textup{i}^{\delta_{a,2}}
\theta_{a+1}(\eta|\tau)}{\theta_1(2 z|\tau) } \cdot
\\  \nonumber &&
\cdot \Bigl[\,\theta_{a+1} \left(2z-n\eta+\eta|\tau\right)e^{\eta\partial_z}
- \theta_{a+1}\left(-2z-n\eta+\eta|\tau\right)e^{-\eta\partial_z}\Bigl]
\left[\theta_3 \left(z | {\textstyle\frac{\tau}{2}}\right)\right]^j \,
\left[\theta_4 \left(z | {\textstyle\frac{\tau}{2}}\right)\right]^{n-1-j},
\nonumber\end{eqnarray}
where $\mu_a=(-1)^{m/2+\ell+1}$ for $a=0,3$ and $\mu_a=(-1)^{\ell+1}$ for $a=1,2$.
It can be seen that the right-hand side contains the standard
Sklyanin algebra generators with $g=n\eta$ acting in the space $\Theta^+_{2n-2}(z|\tau)$,
and we obtain the invariance of the space of functions $\varphi_{j,\ell}^{(n,m)}(z)$
under the action of the generators $\mathbf{S}^a(n\eta+m\tau/2)$. Permuting $\eta$ and
$\tau/2$ and also $n$ and $m$, we also obtain the same statement for the generators
$\mathbf{\tilde S}^a(n\eta+m\tau/2)$.

Now let the integer $m$ be odd. Using the quasiperiodicity of the
functions $\theta_{a}$ under shifts of the argument by $\tau(m+1)/2$, we then obtain
\begin{eqnarray} \nonumber && \makebox[-2em]{}
\mathbf{S}^a(n\eta+m\tau/2)\varphi_{j,\ell}^{(n,m)}(z)=\mu_a\,
e^{-\pi\textup{i}(mn+n-1)\eta-\pi\textup{i}(m+1)\tau/2-\pi\textup{i}(m+1)^2\tau/4}
\cdot \\ \nonumber && \cdot
\left[\theta_3 \left(z | \eta \right)\right]^\ell \,
\left[\theta_4 \left(z | \eta \right)\right]^{m-1-\ell}
\frac{\textup{i}^{\delta_{a,2}}
\theta_{a+1}(\eta|\tau)}{\theta_1(2 z|\tau) }
\Bigl[\,e^{2\pi\textup{i}z}\theta_{a+1} \left(2z-n\eta+\eta
+\textstyle{\frac{\tau}{2}} |\tau\right)e^{\eta\partial_z} -
\\  \nonumber &&
- e^{-2\pi\textup{i}z}\theta_{a+1}\left(-2z-n\eta+\eta
+\textstyle{\frac{\tau}{2}}|\tau\right)e^{-\eta\partial_z}\Bigl]
\left[\theta_3 \left(z | {\textstyle\frac{\tau}{2}}\right)\right]^j \,
\left[\theta_4 \left(z | {\textstyle\frac{\tau}{2}}\right)\right]^{n-1-j},
\nonumber\end{eqnarray}
where $\mu_a=(-1)^{(m+1)/2+\ell}$ for $a=0,3$ and $\mu_a=(-1)^\ell$ for $a=1,2$.
Because
$$
e^{\pi\textup{i}z}\theta_{1,2,3,4}(z+\textstyle{\frac{\tau}{2}}|\tau)
\propto\theta_{4,3,2,1}(z|\tau),
$$
we have a permutation of the action of the standard Sklyanin algebra generators
$\mathbf{S}^{0,1,2,3} \to \mathbf{S}^{2,3,1,0}$
in the right-hand side up to some constant factors, and we obtain the
invariance of the space of functions  $\varphi_{j,\ell}^{(n,m)}(z)$ for odd $m$.

For $n,m>0$, we have thus described a finite-dimensional (more precisely, $nm$-dimensional)
invariant subspace for the Sklyanin algebra, which was partially characterized in
\cite{DS}. In it, we also realize a finite-dimensional representation
of the elliptic modular double. For $m=0$ or $n=0$, we have
a finite-dimensional ($n$- or $m$-dimensional) representation for only
one of the Sklyanin subalgebras in the double. The observation that
the Sklyanin algebra has finite-dimensional representation not only for
the spin values $g=n\eta$ but also for the integer lattices
$n\eta+m{\textstyle\frac{\tau}{2}}$ and
$n\eta + m {\textstyle\frac{\tau}{2}}+{\textstyle\frac{1}{2}}$
with $n,m>0$ is a more or less obvious consequence of the
modular doubling suggested in \cite{AA2008} because there exists an
involution permuting Sklyanin subalgebras. This fact was also noted in  \cite{RR}.

We stress that the elliptic modular double uniquely fixes the constructed
finite-dimensional space. If we multiply the functions indicated above
by a periodic function $\psi(z)$ with the period $\eta$ or $\tau/2$,
then the condition of invariance under the action of double generators
and the requirement of the analyticity in the variable $e^{2\pi\textup{i}z}$
lead to the constraints $\psi(z+\eta)=\psi(z+\tau/2)=\psi(z+1)=\psi(z)$.
Because $2\eta, \tau$, and $1$ are incommensurable, we obtain $\psi(z)=const.$

\section{A nonfactored form of the intertwining operator}

We now transform general expression \eqref{genform} to a ``normal ordered" form.
Using relations \eqref{mper}, we can easily show that
$$
\mathrm{A}_k(g + m {\textstyle\frac{\tau}{2}}) = (-1)^{m\delta_{k,4}}\,
e^{-{\textstyle\frac{\pi \textup{i} \tau}{2}}m^2 -2\pi \textup{i} m(g+\eta) + \pi \textup{i} m \eta}\, \cdot e^{\frac{\pi \textup{i} m}{\eta} z^2}\,
\mathrm{A}_k(g)\,e^{-\frac{\pi \textup{i} m}{\eta} z^2}\,,
$$
and therefore
$$
\mathrm{A}_k(n\eta-\eta+ m {\textstyle\frac{\tau}{2}})\cdots
\mathrm{A}_k(\eta+m {\textstyle\frac{\tau}{2}})
\mathrm{A}_k(m {\textstyle\frac{\tau}{2}}) =
$$
$$
= (-1)^{m n \delta_{k,4}}\,
e^{-{\textstyle\frac{\pi \textup{i} \tau}{2}}\,m^2 n - \pi \textup{i} \eta\, n^2 m}
\, \cdot e^{\frac{\pi \textup{i} m}{\eta} z^2}\,
\mathrm{A}_k(n\eta-\eta)\cdots
\mathrm{A}_k(\eta)
\mathrm{A}_k(0)\,e^{-\frac{\pi \textup{i} m}{\eta} z^2}\,.
$$
We can now single out the operators $\mathrm{M}\left(n\eta\right)$ and
$\mathrm{M}\left(m {\textstyle\frac{\tau}{2}}\right)$ in
$\mathrm{M}\left(n\eta + m {\textstyle\frac{\tau}{2}}\right)$,
\begin{eqnarray*} &&
\mathrm{M}\left(n\eta + m {\textstyle\frac{\tau}{2}}\right) = (-1)^{m n \delta_{k,4}}\,
e^{-{\textstyle\frac{\pi \textup{i} \tau}{2}}\,m^2 n - \pi \textup{i} \eta\, n^2 m}
\, \cdot
\\ && \makebox[4em]{}
\cdot e^{\frac{\pi \textup{i} m}{\eta} z^2}\,
\mathrm{M}(n\eta)\,e^{-\frac{\pi \textup{i} m}{\eta} z^2}\,\cdot
\theta_k^{n} \left(z | {\textstyle\frac{\tau}{2}}\right)
\mathrm{M}\left(m {\textstyle\frac{\tau}{2}}\right)
\theta_k^{-n} \left(z | {\textstyle\frac{\tau}{2}}\right)\,.
\end{eqnarray*}
Using the quasiperiodicity properties of the functions
$\theta_k^{-n} \left(z | {\textstyle\frac{\tau}{2}}\right)$,
we can show that
\begin{equation}
\theta_k^{n} \left(z | {\textstyle\frac{\tau}{2}}\right)
\mathrm{M}\left(m {\textstyle\frac{\tau}{2}}\right)
\theta_k^{-n} \left(z | {\textstyle\frac{\tau}{2}}\right)
=(-1)^{mn\delta_{k,4}}e^{-2\pi i nz^2/\tau}
\mathrm{M}\left(m {\textstyle\frac{\tau}{2}}\right)e^{2\pi i nz^2/\tau}.
\label{Mconj}\end{equation}
Substituting this equality in the preceding relation, we obtain
\begin{eqnarray}\label{M=MM}  &&
\mathrm{M}\left(n\eta + m {\textstyle\frac{\tau}{2}}\right) =
e^{-{\textstyle\frac{\pi \textup{i} \tau}{2}}\,m^2 n - \pi \textup{i} \eta\, n^2 m}
e^{\frac{\pi \textup{i} m}{\eta} z^2}\,
\mathrm{M}(n\eta)\,e^{-\pi\textup{i}m\frac{z^2}{\eta}-2\pi\textup{i} n\frac{z^2}{\tau}}
\mathrm{M}\left(m {\textstyle\frac{\tau}{2}}\right)
e^{2\pi\textup{i} n\frac{z^2}{\tau}}=
\\  \nonumber  && \makebox[5.8em]{}
=e^{-{\textstyle\frac{\pi \textup{i} \tau}{2}}\,m^2 n - \pi \textup{i} \eta\, n^2 m}
e^{\frac{2\pi \textup{i} n}{\tau} z^2}\,
\mathrm{M}(m{\textstyle\frac{\tau}{2}})\,e^{-\pi\textup{i}m\frac{z^2}{\eta}-2\pi\textup{i} n\frac{z^2}{\tau}}
\mathrm{M}\left(n\eta\right)e^{\pi\textup{i} m\frac{z^2}{\eta}},
\end{eqnarray}
where the second relation is obtained from the first
by a simple permutation of $n$ and $\eta$ with $m$ and $\tau/2$.

Using the recurrence relations, we can now easily prove the equality \eqref{Mphi}.
For this, we consider the action of $\mathrm{M}\left(m {\textstyle\frac{\tau}{2}}\right)
e^{2\pi\textup{i} n\frac{z^2}{\tau}}$, the far right factor
in $\mathrm{M}\left(n\eta + m {\textstyle\frac{\tau}{2}}\right)$ given by
\eqref{M=MM}, on $\varphi_{j,\ell}^{(n,m)}(z)$,
\begin{eqnarray*} &&
\mathrm{M}\left(m {\textstyle\frac{\tau}{2}}\right)
e^{2\pi\textup{i} n\frac{z^2}{\tau}}\varphi_{j,\ell}^{(n,m)}(z)
=\mathrm{B}_3((m-1){\textstyle\frac{\tau}{2}})\cdots
\mathrm{B}_3((m-\ell){\textstyle\frac{\tau}{2}})
\cdot \\ && \makebox[2em]{} \cdot
\mathrm{B}_4((m-\ell-1)
{\textstyle\frac{\tau}{2}})\cdots \mathrm{B}_4({\textstyle\frac{\tau}{2}})
\mathrm{M}({\textstyle\frac{\tau}{2}}) \,
e^{2\pi\textup{i} n\frac{z^2}{\tau}}\left[\theta_3 \left(z | {\textstyle\frac{\tau}{2}}\right)\right]^j \,
\left[\theta_4 \left(z | {\textstyle\frac{\tau}{2}}\right)\right]^{n-1-j},
\end{eqnarray*}
where we use the second recurrence relation in \eqref{RR}.
Using the explicit form of the $\mathrm{M}({\textstyle\frac{\tau}{2}})$,
arising in the far right position, we obtain
\begin{eqnarray*} &&
\mathrm{M}({\textstyle\frac{\tau}{2}}) \,
e^{2\pi\textup{i} n\frac{z^2}{\tau}}\left[\theta_3 \left(z | {\textstyle\frac{\tau}{2}}\right)\right]^j \,
\left[\theta_4 \left(z | {\textstyle\frac{\tau}{2}}\right)\right]^{n-1-j}
=\frac{c_B}{\theta_1(2z|2\eta)}
e^{2\pi \textup{i} \frac{z^2}{\tau}}
\cdot \\ && \makebox[4em]{} \cdot
\left(e^{\frac{\tau}{2}\partial_z}-
e^{-\frac{\tau}{2}\partial_z}\right)e^{2\pi \textup{i} (n-1)\frac{z^2}{\tau}}
\left[\theta_3 \left(z | {\textstyle\frac{\tau}{2}}\right)\right]^j \,
\left[\theta_4 \left(z | {\textstyle\frac{\tau}{2}}\right)\right]^{n-1-j}=0,
\end{eqnarray*}
i.e., indeed the functions $\varphi_{j,\ell}^{(n,m)}(z)$
define zero modes of the operator
$\mathrm{M}\left(n\eta + m {\textstyle\frac{\tau}{2}}\right)$.
Moreover, we can multiply $\varphi_{j,\ell}^{(n,m)}(z)$ by arbitrary periodic
functions $\psi(z+2\eta)=\psi(z)$ or $\psi(z+\tau)=\psi(z)$ and still have zero modes of
$\mathrm{M}\left(n\eta + m {\textstyle\frac{\tau}{2}}\right)$.

We consider the effect of similarity transformations for both operators
in \eqref{M=MM} in the nonfactored form. The necessary formulas are
\begin{eqnarray*} &&
e^{\frac{\pi \textup{i} (m+1)}{\eta} z^2}\,
e^{(n-2k)\eta \partial_z}\,e^{-\frac{\pi \textup{i} (m+1)}{\eta} z^2}
= e^{-2\pi \textup{i} z\,(m+1)(n-2k)}\,e^{-\pi \textup{i} \eta\,(m+1)(n-2k)^2}\,
e^{(n-2k)\eta \partial_z}\,,
\\ &&
e^{-\frac{2\pi \textup{i} (n-1)}{\tau} z^2}\,
e^{(m-2\ell)\frac{\tau}{2} \partial_z}\,e^{\frac{2\pi \textup{i} (n-1)}{\tau} z^2}
= e^{2\pi \textup{i} z\,(n-1)(m-2\ell)}\,e^{\pi \textup{i} \frac{\tau}{2}\,(n-1)(m-2\ell)^2}\,
e^{(m-2\ell)\frac{\tau}{2} \partial_z}\, .
\end{eqnarray*}
We can now write
$$
\mathrm{M}\left(n\eta + m {\textstyle\frac{\tau}{2}}\right) =
e^{-{\textstyle\frac{\pi \textup{i} \tau}{2}}\,m^2 n - \pi \textup{i} \eta\, n^2 m}
\, \cdot\,
\overline{\mathrm{M}}(n\eta)\,\cdot\,
\overline{\mathrm{M}}'\left(m {\textstyle\frac{\tau}{2}}\right)\,,
$$
where the transformed operators are
$$
\overline{\mathrm{M}}(n \eta) =
\sum_{j=0}^{n}
(-1)^j\, \alpha_j^{(n)}(z)
\,e^{-2\pi \textup{i} z\,(m+1)(n-2j)}\,e^{-\pi \textup{i} \eta\,(m+1)(n-2j)^2}\,e^{(n-2j)\eta \partial_z}\,,
$$
$$
\overline{\mathrm{M}}'\left( m {\textstyle\frac{\tau}{2}}\right) =
\sum_{\ell=0}^{m}
(-1)^\ell\, \mathrm{\beta}^{(m)}_{\ell}(z)\,
e^{2\pi \textup{i} z\,(n-1)(m-2\ell)}\,
e^{\frac{\pi \textup{i} \tau}{2}\,(n-1)(m-2\ell)^2}\,
e^{(m-2\ell)\,\frac{\tau}{2} \partial_z}\,.
$$
We move all shift operators to the right and obtain
\begin{eqnarray*} &&
\mathrm{M}\left(n\eta + m {\textstyle\frac{\tau}{2}}\right) =
e^{-{\textstyle\frac{\pi \textup{i} \tau}{2}}\,m^2 n - \pi \textup{i} \eta\, n^2 m}
\, \cdot\,\sum_{k=0}^{n}
(-1)^k\, \alpha_k^{(n)}(z)\,\,e^{-2\pi \textup{i} z\,(m+1)(n-2k)}\, \cdot
\\ && \makebox[4em]{}
\cdot\, e^{-\pi \textup{i} \eta\,(m+1)(n-2k)^2}\,\cdot
\sum_{\ell=0}^{m}(-1)^\ell\, \mathrm{\beta}^{(m)}_{\ell}\left(z+\eta(n-2k)\right)
\,\cdot\,
\\ && \makebox[4em]{}
\cdot\,
e^{2\pi \textup{i} (z+\eta(n-2k))\,(n-1)(m-2\ell)}\,\cdot
\, e^{\frac{\pi \textup{i} \tau}{2}\,(n-1)(m-2\ell)^2}\,
\,e^{(n-2k)\eta \partial_z}\,e^{(m-2\ell)\,\frac{\tau}{2} \partial_z}\,.
\end{eqnarray*}
Using formulas \eqref{mper} and the explicit expression for the
coefficients $\mathrm{\beta}^{(m)}$, we obtain
$$
\mathrm{\beta}^{(m)}_{\ell}\left(z+\eta(n-2k)\right) =
(-1)^{nm}\,e^{4\pi \textup{i} z\,m(n-2k)}\,
e^{2\pi \textup{i} \eta\,m(n-2k)^2}\,
e^{\pi \textup{i}\tau\,(m-1)(m-2\ell)(n-2k)}\,
\mathrm{\beta}^{(m)}_{\ell}\left(z\right)\,.
$$
Collecting all together, we obtain our operator in the normal ordered form:
\begin{eqnarray}\nonumber  &&
\mathrm{M}\left(n\eta + m {\textstyle\frac{\tau}{2}}\right) =
(-1)^{n m}\, e^{-{\textstyle\frac{\pi \textup{i} \tau}{2}}\,m^2 n - \pi \textup{i} \eta\, n^2 m}
\, \cdot\,\sum_{k=0}^{n}
(-1)^k\, \alpha_k^{(n)}(z)\,\cdot\,
\\ \nonumber && \makebox[4em]{}
\,\cdot
\sum_{\ell=0}^{m}(-1)^\ell\, \mathrm{\beta}^{(m)}_{\ell}\left(z\right)
\,e^{(n-1)(m-2\ell)\left[\frac{\pi \textup{i} \tau}{2}\,(m-2\ell) +
2\pi \textup{i} \left(z+(n-2k)\eta\right)\,\right]}\,\cdot
\\ && \makebox[4em]{}
\cdot\,e^{(m-1)(n-2k)\left[\pi \textup{i} \eta\,(n-2k)
+2\pi \textup{i} \left(z +(m-2\ell){\textstyle\frac{\tau}{2}}\right)\,\right]}
\cdot
e^{\left[(n-2k)\eta +(m-2\ell)\frac{\tau}{2}\right] \partial_z}\,,
\label{Snm} \end{eqnarray}
where we arrange the phase factors in the form resembling the
transformation laws for theta functions of a definite level
under a shift by the period.

For the lattice $g=n\eta+m{\textstyle\frac{\tau}{2}}+{\textstyle\frac{1}{2}}$,
we have
$$
\mathrm{M}\left(n\eta + m{\textstyle\frac{\tau}{2}}+ {\textstyle\frac{1}{2}}\right) =
\mathrm{M}\left(n\eta + m {\textstyle\frac{\tau}{2}}\right)P,
$$
because the operators $\mathrm{M}^{(k)}$  in relation \eqref{inter1/2}
are independent of $k$.

\section{A complete factorization of the intertwining operator on theta functions}

We apply the derived intertwined operators to the functions with special
transformation properties under the shifts on periods $\tau$ and $2\eta$.
We consider the product $F_N(z)\,G_M(z)$, where the functions $F_N(z)$ and
$G_M(z)$ are transformed as (here $a,b\in \mathbb{Z}$)
$$
F_N(z+\tau a) =
e^{-2N a\,\left[\pi \textup{i} \tau\, a + 2\pi \textup{i} z\,\right]}\,F_N(z)\,,
\quad
G_M(z+2\eta b) =
e^{-4 M b\,\left[\pi \textup{i} \eta\, b + \pi \textup{i} z\,\right]}\,G_M(z)\,.
$$
We write all phase factors in a form close to the phases in (\ref{Snm}).
For holomorphic functions, $F_N(z)$ are
theta functions of modulus $\tau$ and have the order $2N$, while $G_M(z)$ are theta
functions of modulus $2\eta$ and have the order $2M$.

We introduce the parameters $\alpha=0,\pm 1$ and $\beta=0,\pm 1$ such that
$n-\alpha$ and $m-\beta$ are always even integers. This means that
$\alpha=0$ for even $n$ and $\alpha=\pm 1$ for odd $n$.
Analogously, $\beta=0$ for even $m$ and $\beta=\pm 1$ for odd $m$.
We can now single out the full period shifts and obtain
\begin{eqnarray*} && \makebox[-2em]{}
F_N\left(z+\beta\,{\textstyle\frac{\tau}{2}} + (n-2k)\eta +(m-\beta-2\ell){\textstyle\frac{\tau}{2}}\,\right)\,
G_M\left(z+\alpha \eta +(n-\alpha-2k)\eta +(m-2\ell){\textstyle\frac{\tau}{2}}\,\right) =
\\  \nonumber && \makebox[2em]{}
= e^{-N (m-\beta-2\ell)\left[\frac{\pi \textup{i} \tau}{2}\,(m-2\ell)
+ 2\pi \textup{i} \left(z+\beta\,{\textstyle\frac{\tau}{2}}
+(n-2k)\eta\right)\,\right]}\,\cdot\,
\\  \nonumber && \makebox[3em]{} \cdot\,
e^{-M (n-\alpha-2k)\left[\pi \textup{i} \eta\,(n-\alpha-2k)
+2\pi \textup{i}\left(z+\alpha\eta+(m-2\ell){\textstyle\frac{\tau}{2}}\right)\,
\right]}\,\cdot
\\  \nonumber && \makebox[3em]{}
\cdot\,
F_N\left(z+\beta\,{\textstyle\frac{\tau}{2}}+(n-2k)\eta\,\right)\,
G_M\left(z+\alpha\eta+(m-2\ell){\textstyle\frac{\tau}{2}}\,\right).
\end{eqnarray*}
Choosing $N=n-1$ and $M=m-1$, we see an almost
complete  cancellation  of the phase factors:
\begin{eqnarray}\label{answ4} &&
\mathrm{M}\left(n\eta + m {\textstyle\frac{\tau}{2}}\right)\,
F_{n-1}(z)\,G_{m-1}(z) =
(-1)^{n m}\, e^{-{\textstyle\frac{\pi \textup{i} \tau}{2}}\,m^2 n - \pi \textup{i} \eta\, n^2 m}
\,\cdot
\\  \nonumber && \makebox[1em]{}
\,\cdot\,e^{\frac{\pi \textup{i} \tau}{2}\beta^2\,(n-1)}
\,e^{2\pi \textup{i} \eta\,\beta\,(n-1)n}
\,\cdot\,e^{2\pi \textup{i} \beta\,(n-1)\,z}\,\,\cdot\,\Bigg[\sum_{k=0}^{n}
(-1)^k\, \alpha_k^{(n)}(z)\,\cdot\,e^{-4\pi \textup{i} \eta\,\beta\,(n-1)k}\,\cdot
\\  \nonumber && \makebox[1em]{}
\,\cdot
F_{n-1}\left(z+\beta\,{\textstyle\frac{\tau}{2}}+(n-2k)\eta\,\right)\Bigg]\,\cdot
\\  \nonumber && \makebox[1em]{}
\cdot\, e^{\pi \textup{i} \eta\alpha^2\,(m-1)}\,e^{\pi \textup{i} \tau\,\alpha\,(m-1)m}
\,\cdot\,e^{2\pi \textup{i} \alpha\,(m-1)\,z}\,\cdot\,\Bigg[\sum_{\ell=0}^{m}(-1)^\ell\, \mathrm{\beta}^{(m)}_{\ell}\left(z\right)
\,\cdot\,e^{-2\pi \textup{i} \tau\,\alpha\,(m-1)\ell}\,\cdot
\\  \nonumber && \makebox[1em]{}
\,\cdot
G_{m-1}\left(z+\alpha\eta+(m-2\ell){\textstyle\frac{\tau}{2}}\,\right)\Bigg]\,.
\end{eqnarray}
The choice $\alpha=\beta=0$ corresponds to even $n$ and $m$, in which case all
phase factors are absent.

We thus see a complete factorization of the intertwining operator, its representation as
a product of two operators acting in different spaces. We note that we have described
the action of $\mathrm{M}$ for arbitrary theta functions and
zero modes of $\mathrm{M}(n\eta + m\tau/2)$, as shown above, single out the spaces
$\Theta_{2n-2}^+(z|\tau)$ for $F_{n-1}$ and $\Theta_{2m-2}^+(z|2\eta)$ for $G_{m-1}$.

We pass to the multiplicative notation, which is more compact
and convenient for analytic reasons \cite{DS}. We first recall that
$$
p = e^{2\pi \textup{i}\tau}, \qquad q = e^{4\pi \textup{i}\eta},
\qquad \theta_1(z|\tau) =\frac{\mathrm{e}^{-\pi \textup{i}z}
\theta\left(e^{2\pi \textup{i} z};p\right)}{\mathrm{R}(\tau)}.
$$
Substituting the last relation in the definition of the elliptic
binomial coefficients, we see that all coefficients $\mathrm{R}(\tau)$
cancel. Simplifying the resulting expression using the relation
$\theta(z; p) = - z\,\theta(z^{-1} ; p)$, we find
$$
\begin{bmatrix}
n \\ k
\end{bmatrix}_{\tau,2\eta} =(-1)^k\,q^{\frac{1}{2}\,k(n+1)}\,\prod_{b=1}^{k}
\frac{\theta\left(q^{b-n-1}\,; p\right)}{\theta\left(q^{b }\,; p\right)}\,.
$$

In the remaining part of the coefficients $\alpha_k^{(n)}(z)$, the
factors $\mathrm{R}(\tau)$ also cancel, and we obtain
$$
c_A^n \cdot
\frac{\theta_1(2 z + 2 \eta(n-2k)\,|\tau)}{\prod_{j = 0}^{n}
\theta_1(2 z - 2\eta(k-j)\, |\tau)} =
q^{\frac{n}{4}}\,q^{\frac{1}{4}(n-2k)(n-1)}\,e^{2\pi \textup{i} \,n z}
\frac{\theta\left(e^{4\pi \textup{i} z}\, q^{n-2k}\,;p\right)}
{\prod_{j = 0}^{n} \theta\left(e^{4\pi \textup{i} z}\, q^{j-k}\,;p\right)}\,.
$$
Transforming the theta functions to an appropriate form
and using the preceding expression for the elliptic binomial coefficients,
we obtain
$$
\alpha_k^{(n)}(z) = \frac{(-1)^{k+1}\,q^{\frac{n^2}{4}+n(k+1)}\,e^{2\pi \textup{i}\,(n+2) z}}
{\prod_{j = 0}^{n} \theta\left(e^{4\pi \textup{i} z}\, q^{j}\,;p\right)}
\cdot
\,\theta\left(e^{-4\pi \textup{i} z}\, q^{2k-n}\,;p\right)\,
\prod_{b=1}^{k}
\frac{\theta\left(e^{-4\pi \textup{i} z}\, q^{b-n-1}\,,q^{b-n-1}\,; p\right)}
{\theta\left(e^{-4\pi \textup{i} z}\, q^{b}\,,q^{b }\,; p\right)}\,.
$$

Finally, we give the explicit form of the intertwining operator in the
multiplicative notation:
\begin{eqnarray}\label{answ4m} &&
\mathrm{M}\left(n\eta + m {\textstyle\frac{\tau}{2}}\right)\,
F_{n-1}(z)\,G_{m-1}(z) =
(-1)^{n m}\,\,e^{2\pi \textup{i}\,z\,\left[n+m+4+\alpha\,(m-1)+\beta\,(n-1)\right]}\,\cdot\,
\\  \nonumber && \makebox[1em]{}
\cdot\,\frac{q^{\frac{(n^2-\alpha^2)(1-m)}{4}+\beta\,\frac{n(n-1)}{2}+n}}
{\prod_{j = 0}^{n} \theta\left(e^{4\pi \textup{i} z}\, q^{j}\,;p\right)}\cdot
\frac{p^{\frac{(m^2-\beta^2)(1-n)}{4}+ \alpha\,\frac{m(m-1)}{2}+m}}
{\prod_{j = 0}^{m} \theta\left(e^{4\pi \textup{i} z}\, p^{j}\,;q\right)}\,\cdot
\\  \nonumber && \makebox[1em]{}
\,\cdot\Bigg[ \sum_{k=0}^{n}
q^{k n\,(1-\beta)+\beta k}\,\theta\left(e^{-4\pi \textup{i} z}\, q^{2k-n}\,;p\right)\,
\prod_{b=1}^{k}
\frac{\theta\left(e^{-4\pi \textup{i} z}\, q^{b-n-1}\,,q^{b-n-1}\,; p\right)}
{\theta\left(e^{-4\pi \textup{i} z}\, q^{b}\,,q^{b }\,; p\right)}\,\cdot
\\  \nonumber && \makebox[4em]{}
\cdot\, F_{n-1}\left(z+\beta\,{\textstyle\frac{\tau}{2}}+(n-2k)\eta\,\right)\Bigg]\cdot\,
\\  \nonumber && \makebox[1em]{}
\,\cdot\Bigg[
\sum_{\ell=0}^{m}
p^{\ell m\,(1-\alpha) + \alpha \ell }\,\theta\left(e^{-4\pi \textup{i} z}\, p^{2\ell-m}\,;q\right)\,
\prod_{b=1}^{\ell}
\frac{\theta\left(e^{-4\pi \textup{i} z}\, p^{b-m-1}\,,p^{b-m-1}\,; q\right)}
{\theta\left(e^{-4\pi \textup{i} z}\, p^{b}\,,p^{b }\,; q\right)}\,\cdot
\\  \nonumber && \makebox[4em]{}
\cdot\, G_{m-1}\left(z+\alpha\eta+(m-2\ell){\textstyle\frac{\tau}{2}}\,\right)\Bigg].
\end{eqnarray}
Choosing $\alpha=\beta=1$ in this formula, we obtain
the intertwining operator derived in \cite{DS} using residue calculus
(after the change of notation $w=e^{-2\pi \textup{i} z}$, $n=2\ell_q+1$, and $m=2\ell_p+1$):
\begin{eqnarray}\nonumber  &&
\mathrm{M}\left(n\eta + m {\textstyle\frac{\tau}{2}}\right)\,
F_{n-1}(z)\,G_{m-1}(z) =
\frac{(-1)^{n m}\,\,e^{4\pi \textup{i}\,z\,\left[n+m+1\right]}
q^{\frac{(n^2-1)(1-m)}{4}+\frac{(n+1)n}{2}}
p^{\frac{(m^2-1)(1-n)}{4}+ \frac{m(m+1)}{2}} }
{\prod_{j = 0}^{n} \theta\left(e^{4\pi \textup{i} z}\, q^{j}\,;p\right)
\prod_{j = 0}^{m} \theta\left(e^{4\pi \textup{i} z}\, p^{j}\,;q\right)}\,\cdot
\\  \label{answ4m11} && \makebox[0em]{}
\cdot\Bigg[\sum_{k=0}^{n}
q^{k}\,\theta\left(e^{-4\pi \textup{i} z}\, q^{2k-n}\,;p\right)\,
\prod_{b=1}^{k}
\frac{\theta\left(e^{-4\pi \textup{i} z}\, q^{b-n-1}\,,q^{b-n-1}\,; p\right)}
{\theta\left(e^{-4\pi \textup{i} z}\, q^{b}\,,q^{b }\,; p\right)}\,
F_{n-1}\left(z+{\textstyle\frac{\tau}{2}}+(n-2k)\eta\,\right)\Bigg]\cdot\,
\\  \nonumber && \makebox[0em]{}
\cdot\Bigg[
\sum_{\ell=0}^{m}
p^{\ell }\,\theta\left(e^{-4\pi \textup{i} z}\, p^{2\ell-m}\,;q\right)\,
\prod_{b=1}^{\ell}
\frac{\theta\left(e^{-4\pi \textup{i} z}\, p^{b-m-1}\,,p^{b-m-1}\,; q\right)}
{\theta\left(e^{-4\pi \textup{i} z}\, p^{b}\,,p^{b }\,; q\right)}\,
G_{m-1}\left(z+\eta+(m-2\ell){\textstyle\frac{\tau}{2}}\,\right)\Bigg].
\end{eqnarray}

As shown above,
$
\mathrm{M}\left(n\eta + m{\textstyle\frac{\tau}{2}}+ {\textstyle\frac{1}{2}}\right) =
\mathrm{M}\left(n\eta + m {\textstyle\frac{\tau}{2}}\right)P,
$
and the action of the intertwining operator in this case therefore has the same form
as given above with the shift of $z$ in the arguments of $F_{n-1}$ and $G_{m-1}$
functions by 1/2.

\section{The intertwining operator on the dual lattice}

In more detail, we consider what happens with the inversion relation
for intertwining operator \eqref{inv} as we approach the lattices
$g=n\eta+m\tau/2$ and $g=\frac{1}{2}+n\eta+m\tau/2, \, n,m\in\Z$.
Intertwining operator \eqref{S1fin}
is not well defined for $g=-n\eta-m\tau/2$
and $g=\frac{1}{2}-n\eta-m\tau/2,$ where $n,m\in\Z_{\geq0}$
except at the points $(n,m)=(0,0),\, (0,1),\, (1,0)$.
Indeed, the factor $1/\Gamma(-2g)$ in it has poles at these lattice
points, while the integral operator part remains well defined
because $|t|<1, \, t=e^{-2\pi \textup{i}g}$. It is therefore convenient
to introduce the renormalized intertwining operator
\begin{equation}
[\mathrm{M}_{ren}(g)\Phi](z)= \frac{(q;q)_\infty\,(p;p)_\infty}{2}
\int_0^1 \frac{\Gamma(\pm z\pm x -g)}
{\Gamma(\pm 2x)}\Phi(x)dx.
\label{M_ren}\end{equation}
Obviously, we still have the intertwining relations
\begin{equation}
\mathrm{M}_{ren}(g)\,\mathbf{S}^a(g) =
\mathbf{S}^a(-g)\, \mathrm{M}_{ren}(g), \qquad
\mathrm{M}_{ren}(g)\,\mathbf{\tilde S}^a(g) =
\mathbf{\tilde S}^a(-g)\,\mathrm{M}_{ren}(g)\,.
\label{inter2}\end{equation}
Contiguous relations \eqref{recrelA} are now modified to the form
\begin{eqnarray*} &&
\mathrm{A}_k(g)\,\mathrm{M}_{ren}(g) =
\theta(e^{-4\pi\textup{i}(g+\eta)};p)\,
\mathrm{M}_{ren}(g+\eta)\,\theta_k (z|\textstyle\frac{\tau}{2})\,,
\\ &&
\mathrm{B}_k(g)\,\mathrm{M}_{ren}(g) = \theta(e^{-4\pi\textup{i}(g+\tau/2)};q)\,
\mathrm{M}_{ren}\left(g+{\textstyle\frac{\tau}{2}}\right)\,
\theta_k \left(z | \eta\right)\,,
\end{eqnarray*}
where the operators $\mathrm{A}_k(g)$  and $\mathrm{B}_k(g)$ have the same form as before.
We can also modify partner relations \eqref{RR2} similarly.

It is easy to derive the explicit expression for the intertwining operator
$\mathrm{M}_{ren}(-n\eta-m\textstyle{\frac{\tau}{2}})$, $\,n,m\in\Z_{\geq 0}$,
$\, (n,m)\neq(0,0)$. In the multiplicative notation, we have
\begin{eqnarray}\nonumber &&
[\mathrm{M}_{ren}(-n\eta-m{\textstyle \frac{\tau}{2}})\Phi](Z)=
\frac{(p;p)_\infty(q;q)_\infty}{4\pi \textup{i}}
\int_{\mathbb{T}}\frac{dX}{X}\frac{\Phi(X)}{\Gamma_{\! p,q}(X^{\pm2})} \cdot\,
\\ && \makebox[2em]{} \cdot\,
\frac{\prod_{i=0}^{n-1}\theta(q^{i-n/2}p^{m/2}XZ^{\pm1};p)
\prod_{k=0}^{m-1}\theta(q^{-n/2}p^{k-m/2}XZ^{\pm1};q)}
{\theta(q^{n/2}p^{m/2}X^{-1}Z^{\pm1};p)\theta(q^{-n/2}p^{-m/2}XZ^{\pm1};q)}.%
\label{Mdual}\end{eqnarray}

We now consider recurrence relation  \eqref{recrelA} for $g=-\eta$,
$$
\mathrm{A}_k(-\eta)\,\mathrm{M}(-\eta)
= \mathrm{M}(0)\,\bar{\theta}_k (z)=\bar{\theta}_k (z)\,,
$$
where
$$
\mathrm{A}_k(-\eta)= e^{\pi \textup{i}\frac{z^2}{ \eta}}\,
\frac{c_A\theta_k \left(z| {\textstyle\frac{\tau}{2}}\right)  }{\theta_1(2z | \tau)}
\left[ e^{\eta \partial_z} - e^{-\eta \partial_z}
\right]\, e^{-\pi \textup{i}\frac{z^2}{ \eta}}\,.
$$
We thus obtain the operator identity
\begin{equation}
\frac{c_Ae^{- \pi \textup{i} \eta}}{\theta_1(2z | \tau)}
\left[ e^{- 2\pi \textup{i} z}e^{\eta \partial_z} - e^{2\pi \textup{i} z}e^{-\eta \partial_z}
\right]\,\mathrm{M}(-\eta)=\II\,, \quad \text{or} \quad \mathrm{M}(\eta)\,\mathrm{M}(-\eta)=\II\,,
\label{inv_eta}\end{equation}
which follows from the representation for $\mathrm{M}(\eta)$ in formula \eqref{meta}.
This relation explicitly forbids the operator $\mathrm{M}(-\eta)$ to have zero modes.
We note that if we assert the equality
$\mathrm{M}(-\eta)\,\mathrm{M}(\eta)=\II$ as a similar consequence of
dual contiguous relation \eqref{RR2}, then this is incorrect.
As shown above, the operator $\mathrm{M}(-\eta)$ has a singular integrand
$1/\theta(e^{2\pi\textup{i}(x-\eta\pm z)};q)$, the shifts  $z\to z\pm \eta$
used to derive \eqref{RR2} bring pole-type singularities to the integration
contour, and the recurrence relation does not hold. And, indeed,
the relation $\mathrm{M}(-\eta)\,\mathrm{M}(\eta)=\II$ cannot hold, because
we know that $\mathrm{M}(\eta)$ has nontrivial zero modes.
Similarly, we obtain the equations
$$
\mathrm{M}({\textstyle\frac{1}{2}}+\eta)\,\mathrm{M}({\textstyle\frac{1}{2}}-\eta)=
\mathrm{M}(\textstyle{\frac{\tau}{2}})\,\mathrm{M}(-\textstyle{\frac{\tau}{2}})=
\mathrm{M}(\textstyle{\frac{1+\tau}{2}})\,\mathrm{M}(\textstyle{\frac{1-\tau}{2}})=\II
$$
without simple relations for the products
$\mathrm{M}(-\textstyle{\frac{\tau}{2}})\,\mathrm{M}(\textstyle{\frac{\tau}{2}})$
and $\mathrm{M}(\textstyle{\frac{1-\tau}{2}})\,\mathrm{M}(\textstyle{\frac{1+\tau}{2}})$.

We take now inversion relation \eqref{inv}, multiply it by
$\Gamma(-2g)$, and take the limit $g\to n\eta+m\textstyle{\frac{\tau}{2}}$ or
$g\to \textstyle{\frac{1}{2}}+n\eta+m\textstyle{\frac{\tau}{2}}$,
where $n,m\in\Z_{\geq 0}$.
As mentioned above, $\Gamma(-2g)$ vanishes in this limit except
at the points $(n,m)=(0,0), (1,0), (0,1)$, for which operator \eqref{Mdual}
has a singular integrand. Excluding them, we formally obtain
the operator identities
\begin{equation}
\mathrm{M}(n\eta+m\textstyle{\frac{\tau}{2}})\,
\mathrm{M}_{ren}(-n\eta-m\textstyle{\frac{\tau}{2}})=
\mathrm{M}(\textstyle{\frac{1}{2}}+n\eta+m\textstyle{\frac{\tau}{2}})\,
\mathrm{M}_{ren}(\textstyle{\frac{1}{2}}-n\eta-m\textstyle{\frac{\tau}{2}})=0\, .
\label{zero1}\end{equation}
Similarly, we can multiply  \eqref{inv} by $\Gamma(2g)$ and take the limit
$g\to -n\eta-m\textstyle{\frac{\tau}{2}}$ or
$g\to \textstyle{\frac{1}{2}}-n\eta-m\textstyle{\frac{\tau}{2}}$,
where $n,m\in\Z_{\geq 0}$. Again, except at $(n,m)=(0,0), (1,0), (0,1)$, we obtain
\begin{equation}
\mathrm{M}_{ren}(-n\eta-m\textstyle{\frac{\tau}{2}})\,
\mathrm{M}(n\eta+m\textstyle{\frac{\tau}{2}})=
\mathrm{M}_{ren}(\textstyle{\frac{1}{2}}-n\eta-m\textstyle{\frac{\tau}{2}})\,
\mathrm{M}(\textstyle{\frac{1}{2}}+n\eta+m\textstyle{\frac{\tau}{2}})=0.
\label{zero2}\end{equation}

We prove equalities \eqref{zero1} using the results obtained earlier.
We consider the part of the integrand of operator  \eqref{Mdual}
depending on the ``external" variable $Z$. Let $n,m>0$. Then theta functions
in the denominator cancel and, using the first relation in
\eqref{33}, we can represent the products of theta functions in the numerator
as linear combinations of the functions
$\varphi_{j,\ell}^{(n,m)}(z)$ fixed in \eqref{varphi}. As already shown,
these functions are zero modes of the operator
$\mathrm{M}(n\eta+m\textstyle{\frac{\tau}{2}})$ for $n,m>0$.

We now set $m=0$. We consider the action of $\mathrm{M}(n\eta)$
on the part of the integrand of $\mathrm{M}_{ren}(-n\eta)$ depending on $Z$:
$$
\mathrm{M}(n\eta)\frac{\prod_{i=0}^{n-1}\theta(q^{i-n/2}XZ^{\pm1};p)}
{\theta(q^{n/2}X^{-1}Z^{\pm1};p)\theta(q^{-n/2}XZ^{\pm1};q)}.
$$
Because $\theta(q^{n/2}X^{-1}Z^{\pm1};p)=q^nX^{-2}\theta(q^{-n/2}XZ^{\pm1};p)$,
we can cancel part of the theta functions. In the numerator,
there consequently remains the product
$$
\prod_{i=1}^{n-1}\theta(q^{i-n/2}XZ^{\pm1};p),
$$
which can be represented as a linear combination of the theta functions
$$
\theta_3^j(z|\textstyle{\frac{\tau}{2}})\theta_4^{n-1-j}(z|\textstyle{\frac{\tau}{2}}),
\quad j=0,\ldots, n-1.
$$
Hence, we must show that
$$
\mathrm{M}(n\eta)\frac{\theta_3^j(z|\textstyle{\frac{\tau}{2}})
\theta_4^{n-1-j}(z|\textstyle{\frac{\tau}{2}})}
{\theta(q^{-n/2}XZ^{\pm1};q)}=0.
$$
Indeed, using recurrence relation \eqref{recrelA}, we obtain
\begin{eqnarray*} &&
\mathrm{M}(n\eta)\frac{\theta_3^j(z|\textstyle{\frac{\tau}{2}})
\theta_4^{n-1-j}(z|\textstyle{\frac{\tau}{2}})}
{\theta(XZ^{\pm1};q)}=
\\ && \makebox[2em]{}
=\mathrm{A}_3((n-1)\eta)\cdots\mathrm{A}_3((n-j)\eta)
\mathrm{A}_4((n-j-1)\eta)\cdots\mathrm{A}_4(\eta)\mathrm{M}(\eta)
\frac{1}{\theta(XZ^{\pm1};q)}.
\end{eqnarray*}
Using the explicit form of the operator $\mathrm{M}(\eta)$, we obtain
\begin{eqnarray*} &&
\mathrm{M}(\eta)\frac{1}{\theta(XZ^{\pm1};q)}=\frac{1}{\theta(X(q^{1/2}Z)^{\pm1};q)}
\frac{c_Ae^{\pi\textup{i}\frac{z^2}{\eta}}}{\theta_1(2z|\tau)} \, \cdot
\\ &&  \makebox[4em]{} \cdot\,
\left(e^{\eta\partial_z}-\frac{\theta(X(q^{1/2}Z)^{\pm1};q)}
{\theta(X(q^{-1/2}Z)^{\pm1};q)}e^{-\eta\partial_z}\right)
e^{-\pi\textup{i}\frac{z^2}{\eta}}=0.
\end{eqnarray*}
Because the parameter $X$ is arbitrary, we obtain the needed equality.
The derived result leads to the relation
$\mathrm{M}(n\eta)\mathrm{M}_{ren}(-n\eta)=0$ for $n>1$.
This is not true for $n=1$: because of the presence of the shifts by $\pm \eta$,
theta functions in the denominator lead to a pole appearing
on the integration contour in $\mathrm{M}_{ren}$, and the interchange of
integration and finite-difference operator action used above is therefore not allowed.
Permuting the parameters $\tau$ and $2\eta$, we obtain
a completely analogous picture for the operator $\mathrm{M}(m\textstyle{\frac{\tau}{2}})$.
Equalities \eqref{zero2}, apparently, should follow
from \eqref{zero1} after a finite-difference ``integration by parts".

The whole space $\mathrm{Ker}\, \mathrm{M}(g)$ is quite big.
In addition to the zero modes of  $\mathrm{M}(n\eta+m\textstyle{\frac{\tau}{2}})$
with $n,m\in\Z_{\geq 0}$, which we investigated, there exist
zero modes of the operator $\mathrm{M}_{ren}(-n\eta-m\textstyle{\frac{\tau}{2}})$.
At the moment, we cannot completely describe  $\mathrm{Ker}\, \mathrm{M}(g)$
in a closed form.
But we see that the finite-dimensional space of zero modes composed of the products
of theta functions, which we have described above,
can be uniquely characterized as the intersection of two invariant subspaces:
$$
\text{products of theta functions} = \mathrm{Ker}\, \mathrm{M}(g)
\cap \mathrm{Im}\,\mathrm{M}_{ren}(-g),
$$
where
$$
g= n\eta+m\textstyle{\frac{\tau}{2}}, \
\textstyle{\frac{1}{2}}+n\eta+m\textstyle{\frac{\tau}{2}},
\ \ n,m\in\Z_{>0}.
$$
Confirmation of this assertion also follows from the explicit form of
operator \eqref{Mdual}. Using formulas \eqref{33}, we can
factor out the dependence on $z$ in the products of even theta functions
in the integrand in \eqref{Mdual}.
Following the discussion of Sec. 5, we can conclude that
$\mathrm{M}_{ren}(-n\eta-m\textstyle{\frac{\tau}{2}})$  with  $n,m>0$
maps all test functions to the products of even theta functions,
which is precisely the space of zero modes of the operator
$\mathrm{M}(n\eta+m\textstyle{\frac{\tau}{2}})$ with  $n,m>0$
that we have described. In addition, such a representation of
operator \eqref{Mdual} should allow describing
the kernel space of $\mathrm{M}_{ren}(-n\eta-m\textstyle{\frac{\tau}{2}})$
itself as the space of functions for which the appropriate set of
``theta-functional moments" vanishes. This observation requires a
further detailed investigation.

\section{Conclusion}

We have described a finite-dimensional space of zero modes of
the integral operator $\mathrm{M}(g)$ arising for two discrete
spin lattices $g$. In particular, we derived an
explicit expression for  $\mathrm{M}(g)$ acting in this space
as a finite-difference operator.  This operator is an intertwiner
for the elliptic modular double, and its zero modes given by products
of theta functions with two different modular parameters define
finite-dimensional representations of this algebra. As a next step, it is necessary
to find a closed form for the action of generators of the Sklyanin algebras
in some basis of this space. Here, we stress that the choice of the basis is a free option,
and it may drastically simplify the situation in some special cases. For instance,
the products of theta functions
\begin{equation}
h_{k}^{(N)}(w;p,q):=\prod_{j=0}^{k-1}\theta(q^jaw^{\pm1};p)\prod_{j=0}^{N-k-1}
\theta(q^jbw^{\pm1};p),\quad  k,N\in\Z_{\geq 0},
\label{dis-bas}\end{equation}
were used as basis vectors of the space $\Theta_{2N}^+$ in
\cite{ros:elementary,ros:sklyanin} for a simplified analysis of
elliptic $6j$-symbols (these functions first arose as
some intertwining vectors in \cite{tak}).
For the elliptic modular double, the two-index basis vectors
$$
h_{kj}^{(N,M)}(w):=h_{k}^{(N)}(w;p,q)h_{j}^{(M)}(w;q,p)
$$
should be considered, and the considerations of \cite{ros:elementary,ros:sklyanin}
should be generalized to such a case.
We expect with an appropriate choice of the measure it should be possible to
derive the two-index biorthogonal functions in \cite{AA2003}.

The identity
$ \mathrm{M}(g) \R_{12}  =\R_{12}' \mathrm{M}(g),$
where $\R_{12}$ is a solution of the YBE derived in \cite{DS} and $\R_{12}'$ is another
similar operator, shows that the kernel space of the operator $\mathrm{M}(g)$
is mapped onto itself by the R-matrix $\R_{12}$.
Therefore, zero modes of  $\mathrm{M}(g)$
form an invariant space for the action of operator $\R_{12}$.
The explicit form of the corresponding finite-dimensional R-matrices
will be considered in a subsequent publication.

\smallskip

{\bf Acknowledgments.} The work of S.E.D. is supported by RFBR
(grants no. 13-01-12405 and 14-01-00341). The work of V.P.S.
is supported by RFBR (grants no. 11-01-00980, 14-01-00474) and
the NRU HSE Scientific fund (grant no 13-09-0133).

\end{document}